\documentclass{amsart}
\usepackage{amssymb}
\usepackage{pstricks}

\newcommand{\mr}[1]{\mathcal{#1}}
\newcommand{\C}{\mathbb{C}}
\newcommand{\R}{\mathbb{R}}
\newcommand{\Z}{\mathbb{Z}}
\newcommand{\F}{\mathbb{F}}

\newcommand{\HH}{\mr{H}}

\newcommand{\Sym}{\mathrm{Sym}}
\newcommand{\Spec}{\mathrm{Spec}}

\newcommand{\Id}{\mathrm{Id}}
\newcommand{\SL}{\mathrm{SL}}
\newcommand{\SO}{\mathrm{SO}}
\newcommand{\St}{\mathrm{St}}

\newtheorem{lem}{Lemma}[section]
\newtheorem{theorem}[lem]{Theorem}
\newtheorem{obs}[lem]{Observation}

\theoremstyle{remark}
\newtheorem{claim}[lem]{Claim}
\newtheorem{cor}[lem]{Corollary}
\newtheorem{remark}[lem]{Remark}
\newtheorem{defi}[lem]{Definition}

\newtheorem{example}[lem]{Example}

\newcommand{\la}{\langle}
\newcommand{\ra}{\rangle}

\newcommand{\an}[2]{\sphericalangle\left(#1,#2\right)}
\newcommand{\anone}[1]{\sphericalangle\left\{#1\right\}}
\newcommand{\e}{\varepsilon}
\newcommand{\kc}{\kappa} 
\DeclareMathOperator{\ssum}{sum}
\newcommand{\bolde}{\pmb{\e}}
\newcommand{\vsp}{\rule{0pt}{2.5ex}}

\begin{document}
\title{Subspace Arrangements and Property T}
\author{Martin Kassabov}
\thanks{The author was supported in part by the
NSF grants DMS~
0635607 and~0900932}
\dedicatory{To Fritz Grunewald on the occasion of his 60th birthday}
\begin{abstract}
We reformulate and extend the geometric method for proving Kazhdan property T developed by
Dymara and Januszkiewicz and used by Ershov and Jaikin. The main result says that a group $G$
generated by finite subgroups $G_i$ has property T
if the group generated by each pair of subgroups
has property T and sufficiently large Kazhdan constant. Essentially, the same result
was proven by Dymara and Januszkiewicz, however our bound for ``sufficiently large'' is significantly better.

As an application of this result,
we give exact bounds for the Kazhdan constants and the spectral gaps of the random walks on
any finite Coxeter group with respect to the standard generating set, which generalizes
a result of
Bacher and de la~Hapre.
\end{abstract}

\maketitle

\section{Introduction}

One of the aims of this paper is to explain the author's interpretation of
the method for proving property T developed by Dymara and Januszkiewicz in~\cite{DJ}.
This method reduces proving 
property T of a group $G$
to ``local representation theory''
and geometry of configurations of subspaces of a Hilbert space. Here by ``local representation theory''
we mean studying the representations of (relatively) small subgroups in the group $G$.
The second part of the program can be reduced to an 
optimization problem in
some finite dimensional space, however in almost all cases the dimension is too big and
this problem can not be approached directly. Instead methods from linear algebra and
graph theory are used 
(see~\cite{EJ,EJK}).
One unfortunate side effect is that
the simple geometric idea behind this approach gets ``hidden'' in the computations.

The main new result in this paper is a solution of the resulting optimization
problem in one (relatively easy) special case. In some sense, our solution is optimal,
which allows us to get exact bounds for the Kazhdan constants and spectral gaps
in several situations (in the case of any finite Coxeter groups with respect to the
standard generating set or the group $\SO(n)$ for some specific generating set).
The majority of the results were known previously (see~\cite{BdH,DS}),
however
we were able to obtain them almost without using any representation theory. 
In particular we do not use the classification of the irreducible representations
of symmetric groups nor any character estimates.


\medskip

Let us recall the definition of Kazhdan property T: a unit vector $v$ in a unitary
representation of a group $G$ is called $\e$-almost invariant under a generating set
$S$ if $\| g(v) - v \| \leq \e$ for any $g \in S$. One 
way to construct
almost invariant vectors is to take a small perturbations of invariant vectors.
A group has property T, if this is essentially the only way to construct almost invariant vectors.
More precisely: 
\begin{defi}
The Kazhdan constant, denoted by $\kc(G,S)$, of a group $G$
with respect to a generating set $S$ is the largest $\e$
such that the existence of an $\e$-almost invariant vector in a unitary representation
implies the existence of an invariant vector. A finitely generated  discrete%
\footnote{This definition can be extended to locally compact groups by replacing
finite generating set with a compact generating set.
}
group $G$ is said
to have Kazhdan property T if for some (equivalently any) finite generating
set $S$ the Kazhdan constant is positive.
\end{defi}
It follows almost immediately from the above definition that any finite group
has property T. However, computing the Kazhdan constants
even for finite groups is very difficult and there are only a few cases where exact values
are known~\cite{BdH}.
It is well known~\cite{BdlHV,Kaz}
that many infinite groups also have property T, for example any lattice in high rank
Lie group has property T---a typical examples are the groups $\SL_n(\Z)$ and $\SL_n(\F_p[t])$
for $n \geq 3$. Usually, this is proved using the representation theory
of the ambient Lie group~\cite{Kaz}, but this approach
does not produce any bounds for the Kazhdan constants of this groups.
In the last 10 years several algebraic methods for proving property T have
been developed~\cite{EJ, KSLn,KN, Sh, ShICM}. One main advantage of these methods is that
they provide explicit bounds for the Kazhdan constants of these groups,
an other 
is that these methods are applicable 
in a more general setting.

\bigskip

\noindent
\begin{minipage}[b]{0.65\linewidth}
\setlength{\parindent}{1em} 
\indent
One of the ``smallest'' groups which does not have property T is the infinite dihedral group
\[
D_\infty \simeq \la a,b \mid a^2=b^2=1 \ra.
\]
The failure of property T can be easily seen
using $2$ dimensional representations of $D_\infty$. Let $l_a$ and $l_b$ be
two different lines in the Euclidean plane $\R^2$. Such two lines define a representation of
$D_\infty$ on $\R^2$, where the generators act as reflection along these lines.
Thus, the  lines $l_a$ and $l_b$ are the fixed subspaces of the subgroups of
$D_\infty$ generated by $a$ and $b$.
\end{minipage}
\psset{unit=.2mm}%
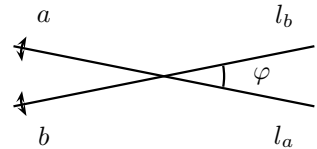
\begin{pspicture}(180,200)(-20,0)
  \psset{linewidth=1.5}
  \psset{linestyle=solid}
  \psline(0,80)(200,120)
  \psline(0,120)(200,80)
  \psarc(100,100){40}{-10}{10}
  \psarc{<->}(25,115){20}{145}{195}
  \psarc{<->}(25,85){20}{165}{215}
  \rput(165,100){$\varphi$}
  \rput(20,140){$a$}
  \rput(20,60){$b$}
  \rput(180,140){$l_b$}
  \rput(180,60){$l_a$}
\end{pspicture}

A quick computation shows that this representation does not contain any invariant vectors
but contains an $\e$-almost invariant vector
where $\e=2\sin \varphi/2$ (here $\varphi$ denotes the angle between the lines $l_a$ and $l_b$). Since the
angle $\varphi$ can be arbitrary small we see that for any $\e >0$ there exists a representation
of $D_\infty$ with $\e$-almost invariant vectors, but without invariant vectors.
In other words,
the group $D_\infty$ does not have property T.

This example suggests 
that a group $G$, generated by two (finite) subgroups
$G_1$ and $G_2$, has property T if and only if there exists a universal bound
for the angle between the fixed spaces $\HH^{G_1}$ and $\HH^{G_2}$ for any
(irreducible) representation $\HH$ of $G$. For example, if we replace $D_\infty$ with $D_n$,
by adding the relation $(ab)^n=1$,
then the reflections along $l_a$ and $l_b$ define representation of $D_n$ if and only if
$\phi = k \pi/n, k\in \Z$ which prevents the angle between $l_a=\HH^{G_1}$ and $l_b=\HH^{G_2}$
from being very small. Of course, this example is not 
useful, since the group
$D_n$ is finite, therefore it has property T.

\medskip

Observation~\ref{ob:main_observation} allows us to
generalize this situation to a group generated by several subgroups, which can be used
to show that some groups have property T by solving a ``geometric optimization'' problem.
Before explaining this reduction and stating the main result in this paper,
we need to look at the angle between two subspaces
from another view point.

We say, that the angle between two closed subspaces $V_1$ and $V_2$ is more than $\varphi$, if for any
vectors $v_i \in V_i$ such that%
\footnote{We need this condition, because we want to measure then angle between subspaces
which have non-trivial intersection.
The motivating example for this definition is the geometric angle
between two planes in a three dimensional
Euclidean space. 
In the theory of Hilbert spaces this angle is known as Friederichs angle~\cite{Deu}.}
each $v_i$ is perpendicular to the intersection $V_1 \cap V_2$,
the angle between $v_1$ and $v_2$ is more then $\varphi$. An equivalent way of
saying this is that for any vector $v\in V$, there is a bound for the distance $d_0(v)$ from $v$ to the
intersection $V_1\cap V_2$ in terms of the distances $d_i(v)$ from $v$ to $V_i$
\[
d_0(v)^2 \leq \frac{1}{1-\cos\varphi} \left( \vsp d_1(v)^2 + d_2(v)^2 \right).
\]

Similar bound can be used to define angle between many
subspaces---we say that the angle $\an{V_1}{V_2,\dots, V_n}$
between the subspaces $V_1,\dots,V_n$ is more than $\varphi$
if for any vector $v$ the square of distance from $v$ to the intersection $\cap V_i$
is bounded by a constant times the sum of the squares of the distances from $v$ to each subspace, i.e.,
\[
d\left(v,\bigcap V_i\right)^2 \leq C_n(\varphi) \sum d(v, V_i)^2,
\]
where $C_n(\varphi)$ is an explicitly defined function.
(See Section~\ref{sec:angles} for a precise definition of the angle between several subspaces.)
Our main result is a sufficient condition when the angle between many subspaces is positive:


\begin{theorem}
\label{th:main-angle-bound}
Let $V_i$ be $n$ closed subspaces in a Hilbert space $\HH$.
If for 
any pair of indexes
$1 \leq i,j \leq n$
we have
$\cos \an{V_i}{V_j} \leq \e_{ij}$ and the
symmetric matrix
\[
A=
\left(
\begin{array}{ccccc}
1 & -\e_{12} & -\e_{13} & \dots & -\e_{1n} \\
-\e_{21} & 1 & -\e_{23} & \dots & -\e_{2n} \\
-\e_{31} & -\e_{32} & 1 & \dots & -\e_{3n} \\
\vdots & \vdots & \vdots & \ddots & \vdots \\
-\e_{n1} & -\e_{n2} & -\e_{n3} & \dots & 1
\end{array}
\right)
\]
is positive definite.
Then, the angle $\an{V_1}{V_2,\dots, V_n} \geq \varphi >0$,
where the constant $\varphi$ depends only on the matrix $A$.

Moreover, if the matrix $A$ is not positive definite then
there exist a Hilbert space $\HH$ and closed subspaces
$V_i$ such that $\cos \an{V_i}{V_j} \leq \e_{ij}$,
but the angle $\an{V_1}{V_2,\dots, V_n}$ is equal to $0$.
\end{theorem}

Weaker forms of Theorem~\ref{th:main-angle-bound} were previously known: Dymara and Januszkiewicz
have proven~\cite{DJ} an analogous statement if $\e_{ij}\leq 12^{-n}$. This result was
improved by Ershov and Jaikin~\cite{EJ} to $\e_{ij}\leq (n-1)^{-1}$. Moreover,
Ershov and Jaikin~\cite[Theorem 5.9]{EJ}
also proved an analog of Theorem~\ref{th:main-angle-bound} in the case $n=3$.

\bigskip

The applications of Theorem~\ref{th:main-angle-bound}
are based on Observation~\ref{ob:main_observation}---a
group $G$ generated by (finite) subgroups $G_i$ has
property T if and only if for any unitary representation of $G$ in $\HH$
there is a bound for the angle between
the subspaces $\HH^{G_i}$, which does not depend on the representation $\HH$.
This
allows one to prove that $G$ has property T
using only information form the representation theory of the groups generated by
$G_i$ and $G_j$. Moreover, a quantitative version of this result
(Theorem~\ref{th:main-distance-estimate}) can be used to obtain good
bounds for the Kazhdan constant and the spectral gap of the Laplacian,
see Theorems~\ref{th:main-app-Coxeter}, \ref{th:main-app-RW-SL_n}
and~\ref{th:main-app-RW-SO_n}---it is remarkable that in some cases
the resulting bounds are sharp.

\begin{theorem}
Let $G$ be a finite Coxeter group with a generating set $S$.
The spectral gap and the Kazhdan constant $\kc(G,S)$ of $G$
can be computed by considering only the defining representation.
In particular the spectral gap of the Laplacian is equal to
$ \frac{4}{n}\sin^2(\pi/2h)$, where $n= |S|$ and $h$ denote the Coxeter number of the group
$G$.
\end{theorem}

This generalizes results by Bacher and de la Harpe~\cite{BdH} and
Bagno~\cite{Bag} and is one of the
few results which provide exact values for the Kazhdan constants
of non-abelian finite groups.
\medskip

Another application of this method is a simplified%
\footnote{A similar proof in the case $p> (n-1)^2$ can be found in~\cite{EJ}.
With a small modification one can extend that proof to the general case.}
proof of the following:

\begin{theorem}
\label{th:main-app-SL_n}
The group
$\SL_n(\F_p[t_1,...,t_k])$ has 
property T, if $p \geq 5$ and $n \geq 3$.
\end{theorem}

The condition $n\geq 3$ is necessary, because the group $\SL_2(\F_p[t])$ does not have property T.
On the other side the condition $p \geq 5$ is redundant---it is even possible to replace $\F_p$ with $\Z$.
However, removing the condition $p \geq 5$ (and replacing $\F_p$ with $\Z$)
requires significant additional work,
see~\cite{EJ} and~\cite{EJK}.
Theorem~\ref{th:main-app-SL_n} also
can be generalized to
show that many Steinberg and Kac-Moody groups have Kazhdan property T.

\bigskip

The proof of Theorem~\ref{th:main-angle-bound}
can also be used to obtain good bounds for the spectral gaps
for some random walks on 
$\SL_n(\F_p)$, $\SO(n)$,
see Theorems~
\ref{th:main-app-RW-SL_n} and~\ref{th:main-app-RW-SO_n},
which in turn can be used to estimate the relaxation and the mixing times of these random walks.
Most of these results are only a slight improvement of previous results~\cite{DS,DS2,KSLn},
however the previous proofs involve completely
different methods and use ``more complicated'' representation theory (at least according to the author).

\subsection*{Notation:}
All the representations in this paper are assumed by unitary.
Throughout the paper $\HH$ will denote an arbitrary a Hilbert space.
As usual
$\la \cdot, \cdot \ra$ will denote the scalar product in $\HH$ and $\| \cdot \|$ will be the norm.
We we use
$\an{v}{w}$ to denote the angle between to nonzero vector in $\HH$.
For a subspace $V$, by $V^{\perp}$ we will denote the orthogonal complement of $V$ in $\HH$
and by $P_V:\HH \to \HH$ the orthogonal projection on $\HH \to V$. The notation
$d_V(v)$ will be used for the
distance between a vector $v\in \HH$ and the subspace $V$, i.e.,
$d_V(v)= \|v - P_v(v)\|$.
We almost never use that $\HH$ is a vector space over the complex number, thus
we often consider $\HH$ only as an Euclidean spaces.
This explains why most of the examples in the paper use
finite dimensional Euclidean spaces (over $\R$)---of course one can ``lift''
all these examples to Hilbert spaces by tensoring with $\C$.

\subsection*{Structure:}
In Section~\ref{sec:how_to_get_T}, we start with an observation, which connects property T and geometry
and use it outline an approach to prove property T for some groups.
The notion of the angle between a collection of subspaces is defined Section~\ref{sec:angles}, which also
contains several properties of this notion.
The following Section~\ref{sec:3_subspaces} contains (relatively easy) technical results about angles between $3$ subspaces and
their intersections. These results are used in Section~\ref{sec:main_result} to prove Theorem~\ref{th:main-angle-bound}. The final Section~\ref{sec:applications} describes
several applications of Theorem~\ref{th:main-angle-bound}.

\section{Strategy for Proving Property T}
\label{sec:how_to_get_T}



A key part of this geometric approach to property T is the following observation,
which allows us to relate property T to 
geometry.

\begin{obs}
\label{ob:main_observation}
Let $G$ be a group and let $G_i$ be a collection of $n$ subgroups in $G$ such that
$G = \la G_1,\dots,G_n \ra$.
Then the Kazhdan constant $\kc(G,\bigcup G_i)$ is strictly positive if and only if
there exists $\alpha >0$ such that
$\an{\HH^{G_1}}{\HH^{G_2},\dots,\HH^{G_n}} > \alpha$ for any unitary representation
$\HH$ of $G$. In particular, if all $G_i$ are finite subgroups and there exists
lower bound for then angle, then $G$ has Kazhdan property T.
\end{obs}
\begin{proof}
Suppose that $\an{\HH^{G_1}}{\HH^{G_2},\dots,\HH^{G_n}} > \alpha$ for any representation $\HH$ of $G$.
Let $\HH$ be arbitrary unitary representation of $G$  and let $v$ be a unit vector in
$\HH$ which is $\e$-almost invariant with respect to $\bigcup G_i$, for some sufficiently small $\e$.
Since each $G_i$ is a subgroup, the almost invariance under $G_i$ implies that the distance from
$v$ to the subspace $\HH^{G_i}$ is less than $\e$. The bound for the angle between the
subspaces $\HH^{G_i}$ implies that
\[
d_{\cap \HH^{G_i}}(v)^2 \leq C \sum d_{\HH^{G_i}}(v)^2 \leq Cn \e^2,
\]
for some constant $C$, which depends only on $\alpha$, but not on the representation $\HH$.
If $\e$ is smaller than $(Cn)^{-1/2}$ then the distance between $v$ and the space of $G$ invariant
vectors $\HH^G = \bigcap \HH^{G_i}$ is less than $1$, therefore there exist nonzero vectors in $\HH^G$.
This shows that
if a representation $\HH$ has an
$\frac{1}{2\sqrt{Cn}}$-almost invariant vector then $\HH$ has an invariant vector, which
is equivalent to $\kc(G,\cup G_i) \geq (Cn)^{-1/2}/2 >0 $.

\smallskip

The other direction is similar---suppose that there is no nontrivial lower bound for the
angle $\an{\HH^{G_1}}{\HH^{G_2},\dots,\HH^{G_n}}$, then for any $C$ there exists a representation $\HH$
and a vector $v$ such that
\[
d_{\cap \HH^{G_i}}(v)^2 \geq C \sum d_{\HH^{G_i}}(v)^2 > 0.
\]
However, the image $\bar v =v / d_{\cap \HH^{G_i}}(v)+ \HH^G$ of vector $v$ in $\HH/\HH^G$
is a unit vector which is
a distance less than $C^{-1/2}$ to each of the subspaces $\HH^{G_i}/\HH^G$. Therefore
$\bar v$ is $2C^{-1/2}$-almost invariant with respect to $\cup G_i$. However, if $C$ is large enough
this contradicts with the $\kc(G,\bigcup G_i) >0$.
\end{proof}

\bigskip

The following outline%
\footnote{This idea goes back to~\cite{DJ} and may be even further to~\cite{Bur,CMS,Ga}.
However, the language used in these papers uses to ``unnecessary geometric objects'',
at least according to the author, which makes these ideas difficult to ``extract''.}
show one possible way to apply the above observation and
use it to prove that some groups $G$ generated by subgroups $G_i$ have property T,
more precisely that the Kazhdan constant $\kc(G,\bigcup G_i)$ is positive.

Briefly the idea is first to extract enough
information (steps 1 and 2) from the ``local representation theory'' of the group $G$
and translate this information to bounds on the angles between some of subspaces $\HH^{G_i}$.
Then one uses ``geometric'' arguments (step 3) to show that these conditions imply a bound for the angle between
all subspaces $\HH^{G_i}$. Although the last step is ``geometric'' in most cases the proof
has 
algebraic flavor and heavily uses linear algebra.

\medskip


The first step in the approach is to study ``local representation theory'' of the group $G$, i.e.,
one can consider the subgroups $G_J = \left\la \bigcup_{j\in J} G_j \right\ra \subset G$ for some subsets
$J \subset \{1,\dots, n\}$. If these groups have property T and there exist good bounds for
the Kazhdan constants $\kc(G_J,\bigcup_{j\in J} G_{j})$, then one can apply
Observation~\ref{ob:main_observation} and translate these into bounds
for the angles
$\anone{\HH^{G_j}\mid
j\in J}$.

The second step (which is optional, but essential for some applications~\cite{EJ,EJK})
is to study the inclusions between the subgroup $G_J$
and translate them into statements about the intersections of the subspaces $\HH^{G_j}$.
For example, the inclusion $G_3 \subset \la G_1, G_2 \ra$ leads to the condition
$\HH^{G_1} \cap \HH^{G_2} \subset \HH^{G_3}$.

The third and final step is to consider all possible configurations of subspaces $V_i$ in
some Hilbert space $\HH$, which satisfy all the conditions found in the first two steps.
If one can prove that
$\an{V_1}{V_2,\dots,V_n} \geq \alpha$ for any ``allowed'' configuration, then we
will get that $\an{\HH^{G_1}}{\HH^{G_2},\dots,\HH^{G_n}} \geq \alpha$ for any representation $\HH$,
which by  Observation~\ref{ob:main_observation} implies that
the group $G$ has a variant of property T. In most cases this geometric problem is
best attacked using tools from linear algebra (and some times graph theory).

Although the third step refers to subspaces in an arbitrary Hilbert space, it is possible to
reduce it to a question about subspaces in a finite dimensional Euclidean space:
the angle $\an{V_1}{V_2,\dots,V_n}$ is determined using the distances between a vector $v$
and the subspaces $V_i$ (of course one need to take a supremum over all vectors $v$).
However, since we work with only one vector at a time, without loss of generality we can assume
that the Hilbert space $\HH$ is spanned by the projections $v_J$ of $v$ onto the subspaces
$\bigcap_{j\in J} V_j$
for all subsets $J\subset \{1,\dots,n \}$, i.e., we can assume that $\dim \HH \leq 2^n$.

Therefore, it is sufficient to consider all possible configuration of subspaces $V_i$ in a $2^n$
dimensional Euclidean space satisfying the conditions found in the first two steps.
Thus, the last step can be reduced
to an optimization problem on some finite dimensional space.
Unfortunately, even for a small $n$, it is very difficult to formulate this optimization
problem and solve it directly---for example proving
Corollary~\ref{cor:angle-bound3} using this idea will involve considering configurations
of 3 three dimensional subspaces in a 6 dimensional Euclidean space. The resulting optimization
problem will involve optimizing function defined on a subset of $7\times 7$ semi-positive definite
symmetric matrices, satisfying certain conditions.
In short, the author does not expect this reduction to be used in practice
(unless it can be implemented on a computer).

\medskip

As we have already mentioned this approach for proving property T is not new and there are several
examples in the literature, where similar program has been carried out: 

Dymara and Januszkiewicz~\cite{DJ} essentially proved that
$\an{V_1}{V_2,\dots,V_n} > \alpha > 0$
if $\cos \an{V_i}{V_j} <\e$ for any $1\leq i,j \leq n$
and a sufficiently small $\e$, i.e., if the subspaces $V_i$ are pairwise almost perpendicular.
They combined this result with bounds coming from the
representation theory of rank 2 groups over finite fields to prove that
$2$-spherical Kac-Moody groups have property T if the defining field is finite and
sufficiently large.

Ershov and Jaikin~\cite{EJ} proved a spectral criterion for property T for groups having a decomposition
as graph of groups. This criterion can be translated in the language outlined above by considering
the fixed subspaces of all vertex and edge groups---the bounds for the codistances at each vertex
are equivalent to bounds for the angles of between the edge spaces. Also, the graph of groups decomposition
imposes restrictions between the intersections of this subspaces.

Ershov and Jaikin also applied this spectra criterion to improve the
Dymara and Januszkiewicz result mentioned above---they
proved that if $\cos \an{V_i}{V_j} < \frac{1}{n-1}$ for any $1\leq i,j\leq n$ then
$\an{V_1}{V_2,\dots,V_n} > \alpha > 0$ and obtained a precise result in the case $n=3$,
which is used to prove a variant of Theorem~\ref{th:main-app-SL_n}.
One of the main results in this paper, Theorem~\ref{th:main-angle-bound}, improves that result and
provides some geometric interpretation.

The main result in~\cite[Theorem~5.5]{EJ} gives bounds
for the Kazhdan constant for groups ``graded by root systems of type $A_2$''
with respect to the union of their root subgroups. Its proof again ``follows''
the general outline described above, but this is not easily seen since
the conditions found in the first two step are 
complicated and
can not be easily translated in a geometric language.
Instead the proof is written in an ``algebraic'' language and heavily uses
linear algebra and graph theory.
This result is generalized in~\cite{EJK} for groups graded by arbitrary root systems.


\section{Angle between subspaces}
\label{sec:angles}

In this section we define the angle between two (and several) closed subspaces of a Hilbert
space. 
We start with a geometric definition,
which is motivated by the standard notion of angles between lines and planes in a
3 dimensional Euclidean space. Then, we find an equivalent definition using
spectrum of certain operators, which later will be used to define ``angle'' between
several subspaces.%
\footnote{It is not clear what is the precise geometric meaning of the angle between
several subspaces. Our definition is closely related to the notion of a codistance used in~\cite{EJ},
see Remark~\ref{rem:angle-codistance}.}

\subsection{Geometric definition}

We start with the usual definition of an angle between a vector and a closed subspace:
\begin{defi}
The angle $\an{v}{V}$
between a closed subspace $V$ and a nonzero vector $v\in \HH$ is defined to
be the angle between $v$ and its projection onto $V$
(or $\pi/2$ is the projection is zero). Equivalently one can use
\[
\| P_V(v) \| = \|v\| \cos \an{v}{V}.
\]
\end{defi}

It is clear that $\an{v}{V} = 0$ if and only if $v\in V$.
Notice, that if we fix the subspace $V$ the function $v \to \an{v}{V}$
is a continuous function defined on $\HH \setminus \{0\}$.

This definition can be extended to angle between 2 subspaces---there are several
natural ways%
\footnote{These extensions are known as Freiderichs and Dixmer angles~\cite{Deu}.}
to do that, which are equivalent if the two subspaces have trivial intersection.
Our approach is to ``ignore'' the intersection by factoring it out, or equivalently
by considering only the orthogonal complement to the intersection.
The definition used in this paper differs from the similar one used in~\cite{DJ,EJ}.


\begin{defi}
\label{def:angle2_geom}
Let $V_1$ and $V_2$ be two closed subspaces in a Hilbert space $\HH$.
If neither of the spaces $V_i$ is contained in the other one,
then the (Friederichs)
angle between $V_1$ and $V_2$ (denoted by $\an{V_1}{V_2}$) is defined
to be the infinimum of the angles between nonzero vectors $v_1$ and $v_2$,
where $v_i \in V_i$ and
$v_i \perp (V_1 \cap V_2)$, i.e.,
\[
\cos \an{V_1}{V_2} =
\sup\left\{ |\la v_1,v_2\ra|\, \mid\, \vsp \|v_i\|=1, v_i \in V_i, v_i \perp (V_1 \cap V_2) \right\}.
\]
\end{defi}
Alternatively, one can define the angle as infimum of the angles between vector and a subspace:
\begin{lem}
\label{lem:2angle}
The angle between $V_1$ and $V_2$ is equal to the infimum of the angles between $V_2$ and
non-zero (or unit) vectors in $V_1$, which are perpendicular to the intersection, i.e.,
\[
\an{V_1}{V_2}
= \inf \left\{ \an{v}{V_2} \mid \vsp v \in V_1, v\perp (V_1 \cap V_2) \right\}.
\]
\end{lem}
\begin{proof}
This follows from the observation that
if $v\perp (V_1 \cap V_2)$ and $v\in V_1$ than $P_{V_2} (v)\perp (V_1 \cap V_2)$.
\end{proof}

\begin{remark}
We need the conditions neither of $V_i$ is a subset of the other one,
because if
$V_2 \subset V_1$ there are no vectors in $V_2$ which are perpendicular
to the intersection of $V_1$ and $V_2$. However, using Lemma~\ref{lem:2angle} one can see
 that the angle ``should be equal'' to $\pi/2$ in the case when $V_2 \subset V_1$ and $V_2 \not = V_1$.
\end{remark}

\begin{cor}
\label{cor:proj-inequalities}
Let $v_1$ be a nonzero vector in $V$ which is perpendicular to the intersection of $V_1 \cap V_2$.
Then:
\begin{itemize}
 \item[a)]
   $\|P_{V_2}(v_1)\| \leq \|v_1\| \cos \an{V_1}{V_2}$;
 \item[b)]
   $\|P_{V_2^\perp}(v_1)\| \geq \|v_1\| \sin \an{V_1}{V_2}$;
 \item[c)]
   $\|P_{V_2}(v_1)\| \leq \|P_{V_2^\perp}(v_1)\| \cot \an{V_1}{V_2} $.
\end{itemize}
\end{cor}

The next result is well known consequence of the compactness of
unit sphere in a finite dimensional Hilbert space.
\begin{lem}
If one of the subspaces $V_1$ and $V_2$ is finite dimensional then there exist vectors
$v_1,v_2 \perp (V_1 \cap V_2)$, $v_i\in V_i$ such that $ \an{V_1}{V_2}=\an{v_1}{v_2}> 0 $.
\end{lem}
\begin{proof}
Assume that $V_1$ is finite dimensional. By Lemma~\ref{lem:2angle}
the angle between $V_1$ and $V_2$ is equal to the infimum of the
function $v\to \an{v}{V_2}$ on the unit sphere in
$V_1 \cap (V_1 \cap V_2)^{\perp}$. This function is a continuous with
compact domain therefore the infimum is achieved.
\end{proof}

\begin{cor}
If at least one of the subspaces $V_i$ is finite dimensional, then the angle $\an{V_1}{V_2}$ is positive.
It is not difficult to construct examples of the infinite dimensional closed subspaces $V_i$ with trivial intersection such that $\an{V_1}{V_2}=0$.
It is easy to see that if $\an{V_1}{V_2}>0$ then $V_1 + V_2$ is a closed subspace.
\end{cor}

\begin{remark}
If the dimensions of $V_1$, $V_2$ and $\HH$ are finite and fixed, then the
angle between the subspaces $V_1$ and $V_2$ can be considered as a function defined
on the subset of the product of two grammarians.
It is important to notice that this function is NOT continuous.
The dimension of the intersection divides the domain into cells and
each cell is open in it closure.
The restriction of the angle to each cell is a continuous function, which tends to zero as
one approaches the boundary of each cell.
\end{remark}

\medskip

The following lemma plays in central role in the rest of the paper---roughly speaking
it allows for exchange all subspace with their orthogonal complements.
As a consequence, we can
replace intersections of subspaces with their sums.
\begin{lem}
\label{lem:perp}
The angle between the orthogonal complements $V_1^\perp$ and $V_2^\perp$
is equal to the angle between the subspaces $V_1$ and $V_2$.
\end{lem}
\begin{proof}
Assume that $\an{V_1}{V_2} >0$.
Let $v_1'$ is a non-zero vector in $V_1^\perp$ which is perpendicular to $V_1^\perp \cap V_2^\perp$.
Since $\left( V_1^\perp \cap V_2^\perp \right)^{\perp} = \overline{V_1 + V_2} = V_1 + V_2$ we have
$v_1' \in V_1 + V_2$, i.e., one can write $v_1' = u_1 + u_2$, where $u_i \in V_i$.
Without loss of generality, we can assume that $u_2 \perp V_1 \cap V_2$.
Let $w'$ be the unique vector of the form $v_1' + \lambda u_2$ which is perpendicular to $u_2$.

\noindent
\begin{minipage}[b]{0.55\linewidth}
\begin{claim} We have the inequality
\[
\|w'\| \leq \|v_1'\| \cos \an{V_1}{V_2}.
\]
\end{claim}
\begin{proof}
The angle between $v_1'$ and $u_2$ is not more than
\[
\pi/2 -  \an{u_2}{V_1} \geq \pi/2 - \an{V_1}{V_2}.
\]
Therefore, the angle between $w'$ and $v_1'$ is at least $\an{V_1}{V_2}$.
\end{proof}

Notice, that by the construction of $u_2$ and $w'$ we have 
$P_{V_2^\perp}(v_1') = P_{V_2^\perp}(w')$.
Therefore,
\end{minipage}
%
\psset{unit=.2mm}%
\begin{pspicture}(180,230)(-20,-20)
  \psset{linewidth=0.6}
  \psset{linestyle=solid}
  \psline(0,100)(200,100)
  \rput(30,80){$V_1$}
  \psline(0,200)(200,0)
  \rput(170,0){$V_2$}
  \psset{linewidth=2.5}
  \psline{->}(100,100)(180,100)
  \rput(185,110){$u_1$}
  \psline{->}(100,100)(20,180)
  \rput(30,190){$u_2$}
  \psset{linestyle=dashed}
  \psline{->}(100,100)(140,140)
  \rput(100,190){$v_1'$}
  \psline{->}(100,100)(100,180)
  \rput(150,150){$w'$}
\end{pspicture}

\[
\|P_{V_2^\perp}(v'_1) \| \leq \|w'\|\leq \|v'_1\| \cos \an{V_1}{V_2},
\]
i.e.,
the angle between $v'$ and the subspaces $V_2^\perp$
is at least $\an{V_1}{V_2}$. This shows that
$
\an{V_1^\perp}{V_2^\perp}
\geq
\an{V_1}{V_2}
$.
The opposite inequality follow form the symmetry, thus 
we have that
$\an{V_1^\perp}{V_2^\perp}
=
\an{V_1}{V_2}$.
\end{proof}

\begin{remark}
This lemma (and its proof) is essentially the same as~\cite[Lemma 2.4]{EJ}.
The statement is somehow cleaner because we do not need to check whether
the intersections of $V_1 \cap V_2$ and $V_1^\perp \cap V_2^\perp$ are trivial or not.
Exactly the same result can be found in~\cite{Deu}.
\end{remark}

\subsection{Spectral Definition}
Before generalizing the notion of angle to several subspaces, we need to
``replace'' the geometric definition with an ``algebraic'' one.

Consider to addition operator
\begin{align*}
\ssum : V_1 \oplus V_2 &\to \HH
\\
\ssum(v_1,v_2) &= v_1 + v_2
\end{align*}
and its ``square'' $\Sigma = \ssum^*\circ \ssum$, where $\ssum^*: \HH \to V_1 \oplus V_2$
denotes the transpose of the operator $\ssum$.
It is easy to see that
\[
\ssum^*(v) = (P_{V_1}(v), P_{V_2}(v))
\quad \mbox{and}\quad
\Sigma(v_1,v_2) = (P_{V_1}(v_1+v_2), P_{V_2}(v_1+v_2))
.
\]
Here $V_1 \oplus V_2$ denotes the external direct sum of $V_1$ and $V_2$, which is
naturally a Hilbert space.

By definition, $\Sigma$ is a positive self adjoin operator.
The triangle inequality implies that the norm of $\Sigma$ is
bounded above by $2$, and the point $2$ appears in the discrete spectrum
of $\Sigma$ if and only if the subspaces $V_1$ and $V_2$
have non-trivial intersection.

\begin{lem}
\label{lem:specsum}
The spectrum of $\Sigma$ is invariant under reflection at $1$, i.e.,
$\lambda \in \Spec(\Sigma)$ if and only if $2-\lambda \in \Spec(\Sigma)$.
\end{lem}
\begin{proof}
A direct computation shows that if $(v_1,v_2)$ is an eigenvector of $\Sigma$
with eigenvalue $\lambda$ then
$(v_1,- v_2)$ is also an eigenvector but with eigenvalue $2-\lambda$.
This shows that the discrete part of the spectrum of $\Sigma$
is invariant under the reflection. Using ``approximate eigenvectors''
one can obtain similar result for the continuous part of the spectrum.
\end{proof}

\begin{lem}
\label{lem:spec2bound}
The spectrum of $\Sigma$ is contained in
\[
\left[\vsp 1-\cos \an{V_1}{V_2}, 1 + \cos \an{V_1}{V_2}\right]
\cup \left\{0,2\right\}.
\]
More over the points $1 \pm \cos\an{V_1}{V_2}$ are in the spectrum.
\end{lem}
\begin{proof}
The eigenspace $W_2$ corresponding to the eigenvalue $2$ consists of all vectors
$(v,v)$ for $v \in V_1 \cap V_2$, similarly the eigenspace $W_0$
corresponding to $0$ consists of all vectors $(v,-v)$ for $v \in V_1 \cap V_2$.
Let $w=(v_1,v_2) \in V_1 \oplus V_2$ be a vector perpendicular to $W_0 \oplus W_2$,
 i.e., $v_1,v_2 \perp V_1 \cap V_2$. Then
\begin{align*}
\la \Sigma(w), w \ra & =
\| v_1 + v_2\|^2 = \|v_1\|^2 + \|v_2\|^2+ 2 \|v_1\|.\|v_2\| \cos (\phi) \leq
\\
& \leq
\left(1+|\cos (\phi)| \right)\left(\vsp \|v_1\|^2 + \|v_2\|^2\right) -
|\cos (\phi)|\left(\vsp \|v_1\| - \|v_2\|\right)^2
\leq
\\
& \leq \left(1+|\cos (\phi)|\right) \| w\| ^2,
\end{align*}
where $\phi$ denotes the angle between the vectors $v_1$ and $v_2$.
Thus, on $(W_0 + W_2)^\perp$ the operator $\Sigma$ is bounded
by $1+\cos \an{V_1}{V_2}$, i.e., the spectrum of $\Sigma$ is
contained in $\left[0,1+\cos \an{V_1}{V_2)} \right] \cup \{2\}$.
By Lemma~\ref{lem:specsum} the interval $\left(0,1-\cos \an{V_1}{V_2}\right)$
is also not part of the spectrum.

Let $v_{1,i}$ and $v_{2,i}$ be sequences of unit vectors,
in $V_1$ and $V_2$ respectively, perpendicular to $V_1\cap V_2$
such that $\an{v_{1,i}}{v_{2,i}} \to \an{V_1}{V_2}$.
By the above computation we see that
\[
\left\la \Sigma(w_i), w_i \right\ra \to
\left(\vsp 1+ \cos \an{V_1}{V_2}\right)\|w_i\|^2,
\]
 where
$w_i=(v_{1,i},v_{2,i})$.
Thus, the point $1+\cos \an{V_1}{V_2}$ is in $\Spec(\Sigma)$.
\end{proof}

\begin{remark}
The above two lemmas,  together with the observation that for any two close subspaces $V_1$  and $V_2$ and
any $v\in \HH$ one has the equality
$$
\|\ssum^*_{V_1,V_2}(v) \|^2 + \|\ssum^*_{V_1^\perp,V_2^\perp}(v) \|^2 = 2 \| v \|^2,
$$
give an alternative proof of Lemma~\ref{lem:perp}.
\end{remark}
Lemma~\ref{lem:spec2bound} allows us to define the angle between
$V_1$ and $V_2$ using the spectral gap of the operator $\Sigma$:

\begin{defi}
\label{def:angle2_spec}
The angle $\an{V_1}{V_2} \in [0,\pi]$ between two subspaces $V_1$ and $V_2$
is defined by
\[
1 + \cos \an{V_1}{V_2} = \sup \left\{\vsp \Spec(\Sigma) \setminus \{2\} \right\}.
\]
\end{defi}

\begin{remark}
This definition does not require that neither of the subspaces is contained in the other one
and allows us to define the angle between two subspaces in these degenerate cases:
\begin{itemize}
\item[a)] if $V_1\subset V_2$ but $V_1 \not = V_2$ then
the spectrum of $\Sigma$ consists $0,1$ and $2$,
therefore $\an{V_1}{V_2} = \pi/2$.
Notice, that this agrees with the Lemma~\ref{lem:2angle};

\item[b)] if $V_1 = V_2$ then
the spectrum of
$\Sigma$ consists $0$ and $2$,
therefore $\an{V_1}{V_2} = \pi$.
\end{itemize}
\end{remark}

\begin{remark}
The following example shows why the angle is not a continuous function on the product of the grassmanians:
Let us consider the angle between a plane and line in a 3 dimensional space.
If the line and the plane are in general position than the angle is equal to
the ``geometrically defined one'', i.e., the angle between the line and its projection on to the plane.
However, if the line lies in the plane the angle is equal to $\pi/2$.

In term of the spectrum of $\Sigma$ we have:
if the line is not inside the plane 
then the spectrum of $\Sigma$ is $\{1-\cos \phi ,1, 1+ \cos \phi \}$, when $\phi$ goes to $0$ the spectrum
becomes $\{0,1,2\}$ and the spectral gap near $2$ suddenly increases from $1-\cos \phi$ to $1$.
\end{remark}

\begin{remark}
Using the spectral definition of the angle between two subspaces it is very easy to see why
the angle $\an{V_1}{V_2}$ is always positive if both $V_1$ and $V_2$ are finite dimensional.
In this case the domain of $\Sigma$ is a finite dimensional vector space, thus, $\Spec(\Sigma)$
is a finite subset of $[0,2]$, therefore it does not contain the interval $(2-\varepsilon, 2)$ for
some $\varepsilon >0$.
\end{remark}

\begin{remark}
Another way to define the angle between $V_1$ and $V_2$ is to use the spectra of the
operator $\Sigma'=P^*\circ P$ where $P: V_1 \to V_2$ is the orthogonal projection from $V_1$ to $V_2$.
In the non-degenerate cases we have
\[
\cos \an{V_1}{V_2} = \sup \left\{\vsp \Spec(\Sigma') \setminus \{1\} \right\}.
\]
The only reason we chose to use the operator $\Sigma$ instead of $\Sigma'$ is to
preserve the symmetry between the two subspaces.

One minor difference is that if one uses the operator $\Sigma'$ to define the angle the ``natural''
extension of angle to the case $V_1=V_2$ will give a different answer $\an{V_1}{V_2} = \pi/2$.
In this paper, we will not deal with this degenerate case, so it does not matter how one defines the angle
if the two subspaces are equal.
\end{remark}

\subsection{Angle between several subspaces}
An analog of Definition~\ref{def:angle2_spec} can be used to define an ``angle'' between
a collection of several subspaces. It is not clear what is the exact geometric meaning
of the angle defined bellow---the definition is equivalent to the one described in the
introduction, see Remark~\ref{rem:angle-distance}. 
The same notion is studied in~\cite{BGM}, there the authors define several ways
to measure the ``angle'' between several subspaces, out definition of angle 
is equivalent to the \emph{Friederichs number} used in~\cite{BGM}. 


\begin{defi}
\label{def:angle_spec}
Let $V_1,V_2,\dots V_n$ be closed subspaces of a Hilbert space $\HH$.
Let
\begin{align*}
\ssum: V_1 \oplus V_2 \oplus \cdots \oplus V_n & \to \HH
\\
\ssum(v_1,v_2,\dots,v_n) & = v_1+ v_2+ \dots +v_n
\end{align*}
denote the addition operator. The ``angle'' $\an{V_1}{V_2,\dots, V_n}$
between the subspaces $\{V_i\}$ is defined by%
\footnote{
The reason of this strange normalization is to make the angle between $0$ and $\pi/2$, unless all
$V_i$ are equal to each other.
}
\[
1 + (n-1) \cos \an{V_1}{V_2,\dots, V_n} =
\sup \left\{\vsp \Spec(\Sigma) \setminus \{ n\} \right\},
\]
where $\Sigma = \ssum^* \circ \ssum$.
\end{defi}
\begin{remark}
The triangle inequality shows that $\Spec(\Sigma) \subset [0,n]$. The point $n$ is in the
discrete part of
$\Spec(\Sigma)$ if and only if the intersection $\cap V_i$ is not trivial.
Similarly, $0$ is in the spectra if and only if there exist vectors $v_i \in V_i$ which are
linearly dependant. It can be shown 
that unless all $V_i$ are the same the
$\Spec(\Sigma)$ contains points in the interval $[1,n)$, i.e.,
$\an{V_1}{V_2,\dots, V_n} \leq \pi/2$.
\end{remark}
\begin{remark}
\label{rem:angle-codistance}
Our definition of angle is closely related to the notion of co-distance used in~\cite{EJ}.
If the intersection $\bigcap V_i$ is trivial, one has
\[
\cos \an{V_1}{V_2,\dots, V_n} = \frac{n \rho(V_1,\dots,V_n) - 1}{n-1}.
\]
\end{remark}
\begin{example}
\label{ex:perp_planes_lines}
If the spaces $V_i$ are pairwise orthogonal and have pairwise trivial intersections,
for example if $V_i$ are the coordinate lines in $n$ dimensional Euclidean space
then
$\an{V_1}{V_2,\dots,V_n} = \pi/2$.

If the orthogonal complements of $V_i$ are pairwise orthogonal and
have pairwise trivial intersections,
for example if $V_i$ are the coordinate hyperplanes in $n$ dimensional Euclidean space, then
$\cos \an{V_1}{V_2,\dots,V_n} = 1- \frac{1}{n-1}$.

This example shows that the analog of Lemma~\ref{lem:perp} does not hold for $n > 2$, i.e.,
it is not true in general that
$\an{V_1}{V_2,\dots,V_n} = \an{V_1^\perp}{V_2^\perp,\dots,V_n^\perp}$.
\end{example}

\subsection{Distance Estimates}

In this section we will assume that $V_1$ and $V_2$ are two closed subspaces in $\HH$ such that
$\an{V_1}{V_2} > 0$. As mentioned before, this condition implies that the subspace
$V_1 + V_2$ is closed.

For a vector $w \in \HH$ with $w_0,w_1, w_2$ and $w_{12}$ we will denote the
projections of $w$ onto the subspaces $V_1 \cap V_2$, $V_1$, $V_2$ and $V_1+V_2$ respectively.
Similarly with $d_0(w), d_1(w), d_2(w), d_{12}(w)$ we will denote the
distances of $w$ to these $4$ subspaces.

\begin{lem}
\label{lem:dist-estimate-weak}
The distances $d_i(w)$ satisfy the inequality
\[
d_0(w)^2 \leq
\frac{1}{\displaystyle 1 - \cos \an{V_1}{V_2}}
\left(d_1(w)^2 + d_2(w)^2\right).
\]
Moreover, if there exist a constant $K> 0$ such that for any vector $w\in \HH$,
one has $d_0(w)^2 \leq K \left(d_1(w)^2 + d_2(w)^2\right)$ then
\[
\cos \an{V_1}{V_2} \leq 1 - K^{-1}.
\]
\end{lem}
\begin{proof}
The vector $w'=w-w_0$ is perpendicular to the intersection
$V_1 \cap V_2$ which is the eigenspace of $\ssum\circ \ssum^*$ corresponding to the eigenvalue $2$.
Since the spectrum of the operators $\ssum\circ \ssum^*$ and $\ssum^*\circ \ssum=\Sigma$ are almost the same
(they might differ only at $0$) and $w'$ is perpendicular to the eigenspace corresponding to $2$,
we have
\[
\la w', \ssum \circ \ssum^* (w') \ra \leq
\left(1 + \cos \an{V_1}{V_2}\right) \|w'\|^2
\]
The left side is equal to
\begin{align*}
\| \ssum^*(w') \|^2 & = \| \left(P_{V_1}(w'), P_{V_2}(w') \right)\|=
\|P_{V_1}(w')\|^2+ \|P_{V_2}(w')\|^2=
\\
& =
2 \|w'\|^2 - (\|w'\|^2 - \|P_{V_1}(w')\|^2) - (\|w'\|^2 - \|P_{V_2}(w')\|^2)=
\\
& =
2d_0(w)^2 - d_1(w)^2 - d_2(w)^2.
\end{align*}
Therefore,
$
\left(1-\cos \an{V_1}{V_2}\right) d_0(w)^2 \leq d_1(w)^2 + d_2(w)^2
$.

The second part follow from the observation that for any $\e>0$ there exits a vector
$w'\perp (V_1 \cap V_2)$
such that $\| \ssum^*(w) \|^2 > (1 + \cos \an{V_1}{V_2} - \e) \|w\|^2$ and by the above computation
for such vector one has
\[
(1 - \cos \an{V_1}{V_2} - \e)d_0(w)^2 \geq d_1(w)^2 + d_2(w)^2.
\mbox{\qedhere}
\]
\end{proof}

\begin{remark}
\label{rem:angle-distance}
The previous proof can be generalized to the case of $n$ subspaces---let
$d_i(w)$ denote the distance between $w$ and $V_i$, and $d_0(w)$ the
distance between $w$ and $\cap V_i$ then we have
\[
 d_0(w)^2 \leq \frac{1}{(n-1)(1 - \cos \an{V_1}{V_2,\dots,V_n})} \sum d_i(w)^2.
\]
Similarly, any bound of the form
$d_0(w)^2 \leq \frac{1}{(n-1)\e} \sum d_i(w)^2$ valid for all vectors $w\in \HH$
implies that
$\cos \an{V_1}{V_2,\dots,V_n} \leq 1 - \e$.

This explains 
why the definition of the angle give in the introduction
is equivalent to the Definition~\ref{def:angle_spec}, 
see~\cite{BGM} for a detailed proof.
\end{remark}

One can obtain a slightly more precise estimate than Lemma~\ref{lem:dist-estimate-weak}
 using the following lemma:
\begin{lem}
\label{lem:dist:proj_intersection}
The distance $\|w_1-w_0\|$ between the projection of $w$
onto $V_1$ and the intersection $V_1 \cap V_2$ is bounded by
\[
\|w_1-w_0\| \leq \frac{ \cos \an{V_1}{V_2} d_1(w) + d_2(w)}{\sin \an{V_1}{V_2}}.
\]
\end{lem}
\noindent
\begin{minipage}[b]{0.70\linewidth}
\setlength{\parindent}{1em} 
\begin{proof}
Let $w'$ be the vector on
$V_1$ such that $w_{12}-w'$ is in $V_2$ and is perpendicular to $V_1 \cap V_2$.
Then $w_1 - w' = P_{V_1}(w_{12}-w')$ and $w_{12} - w_1 =P_{V_1^\perp}(w_{12}-w')$, thus
by Corollary~\ref{cor:proj-inequalities} we have
\[
\|w_1-w'\| \leq \cot \an{V_1}{V_2} \|w_{12} - w_1\|.
\]
Similarly $w_{12} -w_2 = P_{V_2}(w'-w_0)$, i.e.,
\[
\|w'-w_0\| \leq \frac{1}{\sin \an{V_1}{V_2}} \|w_{12} - w_2\|.
\]
Therefore,
\begin{align*}
&
\|w_1-w_0\| \leq \|w_1-w'\| + \|w'-w_0\| \leq
\\
&
\leq
\cot \an{V_1}{V_2} \|w_{12} - w_1\| +
\frac{1}{\sin \an{V_1}{V_2}} \|w_{12} - w_2\|
\leq
\\
&
\leq
\cot \an{V_1}{V_2} d_1(w) +
\frac{1}{\sin \an{V_1}{V_2}} d_2(w).
\mbox{\qedhere}
\end{align*}
\end{proof}
\end{minipage}
%
\psset{unit=.2mm}%
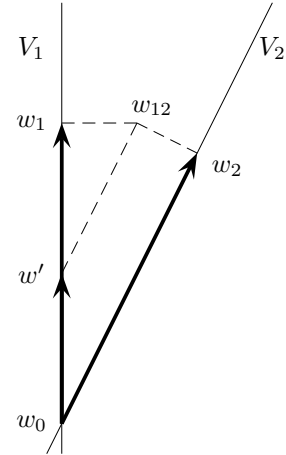
\begin{pspicture}(165,300)(-20,-20)
  \psset{linewidth=0.6}
  \psset{linestyle=solid}
  \psline(20,0)(20,300)
  \rput(0,270){$V_1$}
  \psline(10,0)(160,300)
  \rput(160,270){$V_2$}
  \psset{linewidth=2.5}
  \psset{arrowsize=12}
  \rput(80,230){$w_{12}$}
  \rput(0,20){$w_0$}
  \psline{->}(20,20)(20,220)
  \rput(0,220){$w_1$}
  \psline{->}(20,20)(20,120)
  \rput(0,120){$w'$}
  \psline{->}(20,20)(110,200)
  \rput(130,190){$w_2$}
  \psset{linewidth=0.6}
  \psset{linestyle=dashed}
  \psline(70,220)(20,220)
  \psline(70,220)(20,120)
  \psline(70,220)(110,200)
\end{pspicture}

The following improvement of Lemma~\ref{lem:dist-estimate-weak}
is a special case of 
Theorem~\ref{th:main-distance-estimate}:
\begin{lem}
\label{lem:dist-estimate}
Let $\e=\cos \an{V_1}{V_2} < 1$. Then:
\begin{itemize}
 \item [a)]
$\displaystyle
d_0(w)^2 \leq
\left( \begin{array}{@{}cc@{}} d_1(w) & d_2(w) \end{array} \right)
\left( \begin{array}{@{}cc@{}}
1 & -\e \\
-\e & 1
\end{array} \right)^{-1}
\left( \begin{array}{@{}c@{}} d_1(w) \\ d_2(w) \end{array} \right);
$
\item[b)]
$\displaystyle
\|w_{12}\|^2 \leq
\left( \begin{array}{@{}cc@{}} \|w_1\| & \|w_2\| \end{array} \right)
\left( \begin{array}{@{}cc@{}}
1 & -\e \\
-\e & 1
\end{array} \right)^{-1}
\left( \begin{array}{@{}c@{}} \|w_1\| \\ \|w_2\| \end{array} \right).
$
\end{itemize}
\end{lem}
\begin{proof}
a)
We have
$d_0(w)^2 = \|w_1-w_0\|^2 + d_1(w)^2$. Lemma~\ref{lem:dist:proj_intersection}
gives us a bound for $\|w_1-w_0\|^2$
\begin{align*}
d_0(w)^2
&\leq
\left(\frac{ \cos \an{V_1}{V_2} d_1(w) + d_2(w)} {\sin \an{V_1}{V_2}} \right)^2
+ d_1(w)^2=
\\
& =
\frac{ 1} {\sin^2 \an{V_1}{V_2}}
\left( \vsp d_1(w)^2 + 2 d_1(w)d_2(w)\cos \an{V_1}{V_2} + d_2(w)^2 \right),
\end{align*}
which is equal to bound in the statement of the lemma.

b)
Follow from part a) applied to the subspaces $V_1^\perp$ and $V_2^\perp$.
\end{proof}

\begin{remark}
Part b) of the previous lemma implies that subspace $V_1 + V_2$ is closed if $\an{V_1}{V_2} > 0$.
The analog of this statement is not true for more than $2$ subspaces---there exit closed
subspaces $V_i$ such that $\an{V_1}{V_2,\dots,V_n} > 0$ but $V_1 + V_2 + \cdots + V_n$ is not closed.
However, Remark~\ref{rem:angle-distance} implies that,
if $\an{V_1}{V_2,\dots,V_n} > 0$ then $V_1^\perp + V_2^\perp + \cdots + V_n^\perp$ is
closed and is equal to $\left( \bigcap V_i \right)^{\perp}$.
\end{remark}

\begin{remark}
The bounds in Lemma~\ref{lem:dist:proj_intersection} and Corollary~\ref{lem:dist-estimate}
does not make sense when $V_1=V_2$ and $\an{V_1}{V_2} = \pi$, because the matrix
$\left(\begin{array}{@{}cc@{}} 1 & 1 \\ 1 & 1 \end{array}\right)$ is not invertible.
However, these bounds are still valid if one
resolves the undefined fraction by taking the limit $\an{V_1}{V_2} \to \pi$. The resulting bounds
\[
\|w_1-w_0\|=0
\quad
d_0(w)^2= d_1(w)^2=d_2(w)^2
\quad
\|w_{12}\|^2 = \|w_{1}\|^2 = \|w_{2}\|^2
\]
hold by trivial geometric arguments.
\end{remark}

\section{Three subspaces}
\label{sec:3_subspaces}
In this section, we study configuration of three subspaces in a Hilbert spaces $\HH$.
Corollary~\ref{cor:angle-bound3} gives bound for the angle between the three subspaces in terms of
the angles between each pair. This result is a special case of Theorem~\ref{th:main-angle-bound}.

\subsection{Bounds for the angles}
Let $V_i$ be three closed subspaces in $\HH$ such that $0< \alpha_{ij} \leq \an{V_i}{V_j}$
and $\e_{ij} = \cos \alpha_{ij}$. The next several lemmas give bounds for the angles between
intersections of the subspaces $V_i$-es. Similar results with 
weaker bounds
can be found in~\cite{DJ,EJ}.
Informally these lemmas show that if $W_i$ are
3 planes in $\R^3$ such that $\alpha_{ij} = \an{W_i}{W_j}$, then angle
between two intersections on $V_i$-es is bounded by the angle between the corresponding
intersection of $W_i$-es.

\begin{lem}
\label{lem:angle_line_and_plane}
The angles $\an{V_1+V_2}{V_3}$ and $\an{V_1 \cap V_2}{V_3}$
satisfy the inequalities:
\begin{itemize}
\item[a)]
$
\displaystyle
\cos^2 \an{V_1+V_2}{V_3} \leq
\frac{\e_{13}^2 + \e_{23}^2 + 2\e_{12}\e_{23}\e_{13}}{1- \e_{12}^2}
;$
\item[b)]
$
\displaystyle
\cos^2 \an{V_1 \cap V_2}{V_3}\leq
\frac{\e_{13}^2 + \e_{23}^2 + 2\e_{12}\e_{23}\e_{13}}{1- \e_{12}^2}
.$
\end{itemize}
\end{lem}
\begin{proof}
a)
Let $v_3 \in V_3$ be a vector perpendicular to the intersection $(V_1+V_2) \cap V_3$.
By Lemma~\ref{lem:dist-estimate} we can bound the length of the projection
$P_{V_1+V_2}(v_3)$ of $v_3$ onto $V_1+V_2$ using the length of the
projections $P_{V_1}(v_3)$ and $P_{V_2}(v_3)$. However $v_3$ is perpendicular to both
$V_1 \cap V_3$ and $V_2 \cap V_3$, therefore
we have
\[
\|P_{V_1}(v_3)\| \leq \e_{13}\|v_3\|
\quad\mbox{and}\quad
\|P_{V_2}(v_3)\| \leq \e_{23}\|v_3\|.
\]
Thus,
\begin{align*}
\|P_{V_1+V_2}(v_3)\|^2
& \leq
\|v_3\|^2 \left( \begin{array}{@{}cc@{}} \e_{13} & \e_{23} \end{array} \right)
\left( \begin{array}{@{}cc@{}}
1 & -\e_{12} \\
-\e_{12} & 1
\end{array} \right)^{-1}
\left( \begin{array}{@{}c@{}} \e_{13} \\ \e_{23} \end{array} \right)
=
\\
& =
\frac{\e_{13}^2 + \e_{23}^2 + 2\e_{12}\e_{23}\e_{13}}{1- \e_{12}^2} \|v_3\|^2.
\end{align*}
By definition, any bound of $\|P_{V_1+V_2}(v_3)\|/\|v_3\|$,
which is independent of the vector $v_3$,
is also a bound for
$\cos \an{V_1 +V_2}{V_3}$.

Part b)
follows from part a) by applying it to the subspaces $V_i^\perp$ and using
Lemma~\ref{lem:perp} several times.
\end{proof}
\begin{lem}
\label{lem:angle_line_and_line}
The angles $\an{V_1 +V_3}{V_2+V_3}$ and $\an{V_1 \cap V_3}{ V_2 \cap V_3}$
satisfy the inequalities:
\begin{itemize}
 \item[a)]
$
\displaystyle
\cos \an{V_1 +V_3}{V_2+V_3} \leq
\frac{\e_{12} + \e_{13}\e_{23}}{\sqrt{1- \e_{13}^2}\sqrt{1- \e_{23}^2}};
$

 \item[b)]
$
\displaystyle
\cos \an{V_1 \cap V_3}{ V_2 \cap V_3} \leq
\frac{\e_{12} + \e_{13}\e_{23}}{\sqrt{1- \e_{13}^2}\sqrt{1- \e_{23}^2}}.
$
\end{itemize}
\end{lem}
\begin{proof}
a) Let $w_1$ and $w_2$ are vectors in $V_1+V_3$ and $V_2+V_3$ which are
perpendicular to the intersection $(V_1+V_3)\cap (V_2+V_3)\supset V_3 + (V_1\cap V_2)$.
Then we can find vectors $v_i\in V_i$ and
$v_i \perp V_1 \cap V_2$ such that $w_i = v_i - P_{V_3}(v_i)$.
Then we have
\[
\|w_i\|^2 = \|v_i\|^2 - \|P_{V_3}(v_i)\|^2 \geq \left(1 - \e_{i3}^2\right) \|v_i^2\|
\]
and
\begin{align*}
\la w_1, w_2 \ra
& =
\la v_1 - P_{V_3}(v_1), v_2 - P_{V_3}(v_2) \ra =
\la v_1,v_2 \ra - \la P_{V_3}(v_1), P_{V_3}(v_2) \ra \leq
\\
& \leq
\e_{12}\|v_1\|\, \|v_2\| + \|P_{V_3}(v_1)\|\, \|P_{V_3}(v_2)\|
\leq \left(\e_{12} + \e_{13}\e_{23}\right)\|v_1\|\, \|v_2\|.
\end{align*}
Therefore
\[
\frac{\la w_1, w_2 \ra}{\|w_1\|\, \|w_2\|}
\leq \frac{\e_{12} + \e_{13}\e_{23}}{\sqrt{1- \e_{13}^2}\sqrt{1- \e_{23}^2}}.
\]

Again, part b) can be obtained by applying part a) to the subspaces $V_i^\perp$ and using
Lemma~\ref{lem:perp}.
\end{proof}

\subsection{Relations with spherical geometry}
All the bounds obtained in the previous subsection are ``sharp'' and
have an easy geometric interpretation. 
Let us start with the observation that these bounds are nontrivial only if
$\e_{ij}$ satisfy the inequality:
\[
\e_{12}^2 + \e_{23}^2+ \e_{13}^2 + 2 \e_{12}\e_{23}\e_{13} < 1.
\]
This condition is equivalent to the positive definiteness of matrix
\[
\left( \begin{array}{@{}ccc@{}}
1 & -\e_{12} & -\e_{13}\\
-\e_{12} & 1 & -\e_{23}\\
-\e_{13} & -\e_{23} & 1
\end{array} \right)
.
\]
which is equivalent to $\alpha_{12} + \alpha_{23} + \alpha_{13} > \pi$.

Notice, that this is equivalent to the existence of
3 unit vectors $w_i\in \R^3$ such that $\la w_i, w_j \ra = -\e_{ij}$
(such configuration of vectors is unique up to isometry of $\R^3$).
These vectors define three lines $W_i = \R w_i$ and three planes $W_i' = W_i^\perp$.

It is a well know fact from the spherical geometry that
the cosine of the angle between the line $W_3$ and the plane $W_1+W_2$
is given by the formula in part a) of Lemma~\ref{lem:angle_line_and_plane}.
Similarly the cosine of angle between the planes $W_1+W_3$ and $W_2+W_3$
is equal to the expression in part a) of Lemma~\ref{lem:angle_line_and_line}.
The analogues formulas in parts b) of these lemmas correspond
to the angles between the intersections constructed starting from the planes $W_i'$.

One way to prove the second fact is the use the Gramm Schmidt process.
Let $w'_1=w_1 + \lambda_1 w_3$ and $w'_2 =w_2 + \lambda_2 w_3 $
be the projections of $w_1$ and $w_2$ onto the plane perpendicular to $w_3$, here
\[
\lambda_1 = - \frac{ \la w_1,w_3 \ra}{\la w_3,w_3 \ra} = \e_{13}
\qquad
\lambda_2 = - \frac{ \la w_2,w_3 \ra}{\la w_3,w_3 \ra} = \e_{23}.
\]
These vectors are in the planes $W_1+W_3$ and $W_2+W_3$ and by
construction are perpendicular to their intersection $(W_1+W_3)\cap(W_2+W_3)=W_3$.
Therefore the angle between these vectors is equal to the angle
between the two planes. A direct computations shows that
\[
\cos \an{w'_1}{w'_2} =
\frac{\la w'_1, w'_2\ra }{\|w'_1\|\,\, \| w'_2\| } =
- \frac{\e_{12} + \e_{13}\e_{23}}{\sqrt{1- \e_{13}^2}\sqrt{1- \e_{23}^2}}.
\]


If $\alpha_{12} + \alpha_{23} + \alpha_{13} \leq \pi$ then the bounds in
Lemmas~\ref{lem:angle_line_and_plane} and~\ref{lem:angle_line_and_line} are trivial.
The reason for that is that it is possible to construct an example of subspaces
$V_i$ in $\R^3$ such that $\an{V_i}{V_j} \leq \alpha_{ij}$ where the angles
$\an{V_1}{V_2 \cap V_3}$ and $\an{V_1 \cap V_3}{V_2 \cap V_3}$ are arbitrary
small.

\subsection{Distance Estimates}

Let $V_i$ be three subspaces in a Hilbert space $\HH$. For a vector $v\in \HH$ let
$d_i(v)$ denote the distance between $v$ and the subspace $V_i$, also
let $d_0(v)$ denote the distance between $v$ and the intersection $\bigcap V_i$.
\begin{lem}
\label{lem:dist-estimate3}
If $\e_{12}^2 + \e_{23}^2+ \e_{13}^2 + 2 \e_{12}\e_{23}\e_{13} < 1$,
then
\[
d_0(v)^2 \leq
\left( \begin{array}{@{}ccc@{}} d_1(v) & d_2(v) & d_3(v) \end{array} \right)
\left( \begin{array}{@{}ccc@{}}
1 & -\e_{12} & -\e_{13}\\
-\e_{12} & 1 & -\e_{23}\\
-\e_{13} & -\e_{23} & 1
\end{array} \right)^{-1}
\left( \begin{array}{@{}c@{}} d_1(v) \\ d_2(v) \\ d_3(v) \end{array} \right).
\]
\end{lem}
\begin{proof}
Let $v_3 = P_{V_3}(v)$ denotes the projection of $v$ on to $V_3$.
By construction we have
$d_0(v)^2 = d_3(v)^2 + d_0(v_3)^2$, thus our goal is to bound the distance $d_0(v_3)$.

Let $W_1 = V_1 \cap V_3$ and $W_3 = V_2 \cap V_3$, notice that $W_1\cap W_2 = \cap V_i$
therefore by Lemma~\ref{lem:dist-estimate} we can bound $d_0(v_3)$ using the distances between
$v_3$ and the subspaces $W_1$ and $W_2$, and the angle $\an{W_1}{W_2}$. Lemma~\ref{lem:dist:proj_intersection} and~\ref{lem:angle_line_and_line} provide us with bounds for all these. Substituting everything we get
\begin{align*}
d_0(v)^2
& =
d_0(v_3)^2 + d_3(v)^2\leq
\\
& \leq
\left( \begin{array}{@{}cc@{}} d_{W_1}(v_3) & d_{W_2}(v_3) \end{array} \right)
\left( \begin{array}{@{}cc@{}}
1 & -\e \\
-\e & 1
\end{array} \right)^{-1}
\left( \begin{array}{@{}c@{}} d_{W_1}(v_3) \\ d_{W_2}(v_3) \end{array} \right)
+
d_3(v)^2,
\end{align*}
where $\e=\cos \an{W_1}{W_2}$
\[
d_{W_i}(v_3) \leq \frac{\cos \an{V_1}{V_3} d_3(v) + d_i(v) }{\sin \an{V_1}{V_3}}
\]
and
\[
\cos \an{W_1}{W_2} \leq \frac{\cos \an{V_1}{V_2} + \cos \an{V_1}{V_3}\cos \an{V_1}{V_3}}
{\sin \an{V_1}{V_3}\sin \an{V_1}{V_3}}.
\]

A long and boring computation shows that the resulting expression is exactly equal to
formula in the statement of the lemma. The proof of Theorem~\ref{th:main-distance-estimate}
follows the same idea and shows how to avoid doing this long computation.
\end{proof}

\begin{remark}
The bound in the above lemma is sharp---let $W_i'$ be the 3 planes in $\R^3$ such that
the angles between them are equal to $\alpha_{ij}$. It is easy to see, that for any $3$
positive numbers $d_i$ there exists vector $v\in \R^3$ such that $d_{W'_i}(v) = d_i$
and $\|v\|^2$ is given by the formula above.
\end{remark}

An immediate application of the above bound is the following corollary, which is
a special case of Theorem~\ref{th:main-angle-bound}:
\begin{cor}
\label{cor:angle-bound3}
If $\e_{12}^2 + \e_{23}^2+ \e_{13}^2 + 2 \e_{12}\e_{23}\e_{13} < 1$,
then
\[
1 - \cos \an{V_1}{V_2,V_3} \geq \frac{\lambda}{2},
\]
where $\lambda$ is the smallest
eigenvalue of the (positive definite) matrix
\[
A_\e=\left( \begin{array}{@{}ccc@{}}
1 & -\e_{12} & -\e_{13}\\
-\e_{12} & 1 & -\e_{23}\\
-\e_{13} & -\e_{23} & 1
\end{array} \right).
\]
In particular there is a (nontrivial) lower bound for
$\an{V_1}{V_2,V_3}$ which depends only on $\e_{ij}$.
\end{cor}

\section{Main Result}
\label{sec:main_result}

\subsection{Proof of Theorem~\ref{th:main-angle-bound}}
We start with a quantitative variant of Theorem~\ref{th:main-angle-bound},
which will be used in Section~\ref{sec:applications} to obtain 
bounds for Kazhdan constants and spectral gaps.

\begin{theorem}
\label{th:main-distance-estimate}
Let $V_i$ be $n$ closed subspaces of in a Hilbert space $\HH$. Suppose that
the $n\times n$ symmetric matrix
\[
A=
\left(
\begin{array}{ccccc}
1 & -\e_{12} & -\e_{13} & \dots & -\e_{1n} \\
-\e_{21} & 1 & -\e_{23} & \dots & -\e_{2n} \\
-\e_{31} & -\e_{32} & 1 & \dots & -\e_{3n} \\
\vdots & \vdots & \vdots & \ddots & \vdots \\
-\e_{n1} & -\e_{n2} & -\e_{n3} & \dots & 1
\end{array}
\right).
\]
where $\e_{ij} = \cos \an{V_i}{V_j}$,
is positive definite.
Then for any $v\in \HH$ we have
\[
d_0(v)^2 \leq \mathbf{d}_v^t A^{-1}\mathbf{d}_v,
\]
where $d_0(v)$ denotes the distance between $v$ and $\bigcap V_i$ and
$\mathbf{d}_v$ is the column vector with entries the distances $d_{V_i}(v)$.
\end{theorem}
\begin{proof}
The proof is by induction on $n$. The base case $n=2$ is Lemma~\ref{lem:dist-estimate}.
The induction step follows the idea of Lemma~\ref{lem:dist-estimate3}---let
$V_i' = V_i \cap V_n$ and let $v' = P_{V_n}(v)$.
Using Lemma~\ref{lem:angle_line_and_line} one can bound the angles between
$V_i'$ and used these bounds the to form $(n-1) \times (n-1)$ matrix $A'$. Also,
the distances between $v'$ and $V'_i$ can be bounded by Lemma~\ref{lem:dist:proj_intersection}
and these bounds can be combined in a vector $\mathbf{d}_v'$. In order to complete the induction step
we need to show that 1) the matrix $A'$ is positive definite and 2) the equality
\begin{equation}
\label{recursion}
\mathbf{d}_v^t A^{-1}\mathbf{d}_v = \mathbf{d}_v'^t A'^{-1}\mathbf{d}_v' + d_{V_n}(v)^2.
\end{equation}

\medskip

The matrix $A$ can be written as the product
\[
A =
\left(
\begin{array}{@{}cc@{}}
\Id & -{\bolde}_n \\
0 & 1
\end{array}
\right)
\left(
\begin{array}{@{}cc@{}}
\tilde A & 0\\
0 & 1
\end{array}
\right)
\left(
\begin{array}{@{}cc@{}}
\Id & 0 \\
-{\bolde}_n^t & 1
\end{array}
\right),
\]
where $\tilde A$ is $(n-1) \times (n-1)$ matrix with diagonal entries $1-\e_{in}^2$ and
the off diagonal entries $-\e_{ij} - \e_{in}\e_{jn}$.
Here ${\bolde}_n$ denotes the column vector with
entries $\e_{in}$.

The decomposition of $A$ as a product implies that
\[
\mathbf{d}_v^t A^{-1}\mathbf{d}_v =
\mathbf{d}_v^t
\left(
\begin{array}{@{}cc@{}}
\Id & 0 \\
{\bolde}_n^t & 1
\end{array}
\right)
\left(
\begin{array}{@{}cc@{}}
\tilde A & 0\\
0 & 1
\end{array}
\right)^{-1}
\left(
\begin{array}{@{}cc@{}}
\Id & {\bolde}_n \\
0 & 1
\end{array}
\right)
\mathbf{d}_v =
\mathbf{\tilde d}_v^t \tilde A^{-1}\mathbf{\tilde d}_v + d_{V_n}(v)^2,
\]
where $\mathbf{\tilde d}$ is the vector defined by
\[
\left(
\begin{array}{@{}c@{}}
\mathbf{\tilde d}_v \\ d_n
\end{array}
\right)
=
\left(
\begin{array}{@{}cc@{}}
\Id & {\bolde}_n \\
0 & 1
\end{array}
\right)
\mathbf{d}_v.
\]
The equality (\ref{recursion}) follow from the observations, which are immediate consequences of
the definitions of $A'$ and $\mathbf{d}'_v$:
\begin{enumerate}
\item[(a)] the matrices $A'$ and $\tilde A$ are related by $A' = D \tilde A D$, where
$D$ is a diagonal matrix with entries $1/\sqrt{1-\e_{in}}$;
\item[(b)] the vectors $\mathbf{d}'_v$ and $\mathbf{\tilde d}_v$ satisfy
$\mathbf{d}'_v = D \mathbf{\tilde d}_v$,
\end{enumerate}
because they imply that
\[
\mathbf{d}_v'^t A'^{-1}\mathbf{d}_v'=
\mathbf{\tilde d}_v'^t DA'^{-1}D\mathbf{\tilde d}_v'=
\mathbf{\tilde d}_v^t \tilde A^{-1}\mathbf{\tilde d}_v.
\]
The first observation also proves that $A'$ is a positive definite matrix, since $\tilde A$ is.
\end{proof}
\begin{remark}
The geometric interpretation of the above theorem and its proof is the following:
let $w_i$ are unit vectors in $\R^n$ such that $\la w_i, w_j \ra = - \e_{ij}$
(such vectors exists since $A$ is positive definite) and let $w$ be a vector such that
$\la w, w_i \ra = d_{P_{V_i}}(v)$. Then $d_0(v) \leq \| w\|$.

The proof is by induction---the induction step uses Gramm-Schmidt process:
one projects the vectors $w_i$ onto the hyperplane perpendicular to $w_n$
and then ``normalize'' the resulting
vectors $w_i'$ to have a unit length. By the Lemma~\ref{lem:angle_line_and_line} the angle
between $w_i'$ and $w_j'$ is a bound for the angle between the intersections
$V_i \cap V_n$ and $V_j \cap V_n$. Similarly by Lemma~\ref{lem:dist:proj_intersection}
the distance between the projection $P_{V_n}(v)$
to the intersection $V_i \cap V_n$ is bounded by $\frac{\la w', w_i' \ra}{\|w_i'\|^2}$, where
$w'$ is the projection of $w$ onto the hyperplane perpendicular to $w_n$.
The induction step is completed by 
\[
d_0(v)^2 = d_0\left(P_{V_n}(v)\right)^2 + d_{V_n}(v)^2 \leq
\|w'\|^2 + \la w ,w_n \ra^2 = \| w\|^2,
\]
where the inequality follows by the induction assumption.
\end{remark}

\begin{remark}
The inequality in Theorem~\ref{th:main-distance-estimate} can be rephrased as follows:\\
the $(n+1) \times (n+1)$ matrix
\[
B =
\left(
\begin{array}{@{}cc@{}}
d_0(v)^2 & \mathbf{d}_v^t \\
\mathbf{d}_v & A
\end{array}
\right)
\]
is not positive definite.
\end{remark}

\begin{proof}[Proof of Theorem~\ref{th:main-angle-bound}]
If the matrix $A$ is positive definite  Theorem~\ref{th:main-distance-estimate} applies.
Therefore, we have
\[
d_0(v)^2 \leq \mathbf{d}_v^t A^{-1}\mathbf{d}_v \leq \frac{1}{\lambda} \sum d_{V_i}(v)^2,
\]
where $\lambda$ is the smallest eigenvalues of the matrix $A$.
By Remark~\ref{rem:angle-distance} this bound implies that
\[
\cos \an{V_1}{V_2,\dots,V_n} \leq 1 - \frac{\lambda}{n-1},
\quad
\mbox{i.e.,}
\quad
\an{V_1}{V_2,\dots,V_n} \geq \alpha,
\]
where $\alpha = \cos^{-1}\left(1 - \frac{\lambda}{n-1}\right)$
which completes the proof, since $\lambda$ and $\alpha$ depend only on $\e_{ij}$.
\end{proof}

\begin{remark}
In the special case when all $\e_{ij}=\e$ are the same we obtain that if $\e \leq \frac{1}{n-1}$
then $\an{V_1}{V_2,\dots,V_n} > \alpha$ where $\cos \alpha = \frac{n-2}{n-1} + \e$, because the
smallest eigenvalue $\lambda$ of the matrix $A$ is equal to $1-(n-1)\e$,
which is equivalent to Corollary~5.3 from~\cite{EJ}.
\end{remark}
\begin{example}
If $V_i$ are pairwise orthogonal subspaces in $\HH$ then the matrix $A$ is equal to the identity
matrix, and by Theorem~\ref{th:main-angle-bound} we have
\[
\cos \an{V_1}{V_2,\dots,V_n} \leq 1- \frac{1}{n-1}.
\]
In fact, by Example~\ref{ex:perp_planes_lines}
there is an equality if $V_i$ are the coordinate hyperplanes in $n$ dimensional Euclidean space.
\end{example}

\subsection{Geometric Interpretation}
Theorem~\ref{th:main-angle-bound} can be rephrased as follows: let $\Delta$ be a spherical
simplex such the internal angle between any two faces $F_i$ and $F_j$ is equal to $\alpha_{ij}$.
Then the angle between any collection of subspaces $V_i$ such that $\an{V_i}{V_j}\geq \alpha_{ij}$ is
bounded by
\[
\an{V_1}{V_2,\dots,V_n} \geq \an{\tilde F_1}{\tilde F_2,\dots, \tilde F_n} > 0,
\]
where $\tilde F_n$ is the affine subspace which contains the face $F_i$.
A slight modification of the proof also gives that similar inequality holds
for the angle between intersections of $V_i$-es.

Theorem~\ref{th:main-distance-estimate} has a similar interpretation: let $p$ be any point
in the interior of the simplistic cone defined by $\Delta$. Then for any $v\in \HH$
such that $d_{V_i}(v) \leq d_{\tilde F_i}(p)$ we have that
\[
d_{\cap V_i}(v) \leq
d_{\cap \tilde F_i}(p) = \| p\|.
\]

\section{Applications}
\label{sec:applications}

\subsection{Kazhdan constants and Spectral gap for Coxeter groups}
Let $G$ be a finite group generated by a symmetric set $S$, i.e., $S=S^{-1}$.
Let $\pi:G \to U(L^2(G))$ denote the regular representation of the group $G$.
The operator
\[
\Delta_S = \frac{1}{|S|}\sum_{s \in S} \left(\vsp \Id - \pi(s)\right) :L^2(G) \to L^2(G)
\]
is called Laplacian%
\footnote{The operator $\Delta$ can be defined even if the group is not finite, however
in this setting there is not direct connection between $\Delta_S$ and a graph Laplacian.
One can even define $\Delta_\mu$ when $G$ as a group and $\mu$ is a measure on $G$.
In these more general situations there is also a connection between the spectral gap
of $\Delta$ (if positive) and the relaxation time of some random walk.}
on $G$. An equivalent way to define this operator is to take
the graph Laplacian of the Cayley graph of the group $G$ with respect to
the generating set $S$. This operator is positive definite and has an eigenvalue $0$
with multiplicity $1$ (the eigenvector is the constant function).
The 
smallest not-trivial eigenvalue $\lambda_S$ of $\Delta_S$ is called the spectral gap
of the Laplacian and is closely related with the relaxation time of the random
walk on Cayley graph. Thus, bound for the spectral gap can be used to estimate the
mixing time of this random walk.

\medskip

A Coxeter group $G$ generated by a set $S=\{s_1,\dots, s_n\}$ is defined by numbers $m_{ij} \geq 2$ and
has presentation
\[
G\simeq \left\la s_i \mid \vsp s_i^2 =1, (s_is_j)^{m_{ij}}=1 \right\ra.
\]
It is known that $G$ has a \emph{defining representation} on a $n$-dimensional
vectors space $V$ where each generator $s_i$ acts as a reflection with respect to
a hyperplane $V_i$. Moreover if $G$ is finite there is a $G$-invariant Euclidean
structure on $V$ and the angle between the hyperplanes $V_i$ and $V_j$ is equal to
$\pi/m_{ij}$.
\begin{theorem}
\label{th:main-app-Coxeter}
Let $G$ be a finite Coxeter group.
Then, the Kazhdan constant $\kc(G,S)$ and the spectral gap of the Laplacian can be
computed using the defining representation of $G$.
\end{theorem}
\begin{proof}
The group $G$ is generated by $n$ subgroups $G_i = \{1,s_i\}$ of order $2$.
From the presentation of $G$ it is clear that the group generated by $G_i$ and $G_j$
is the dihedral group $D_{m_{ij}}$.
Therefore, 
for any unitary
representation of $G$ in $\HH$ the angle between $\HH^{G_i}$ and $\HH^{G_j}$ is bounded bellow by
$\pi/m_{ij}$. 
It is classical fact~\cite{Hymp}
that the faintness of the group $G$ is equivalent to the
positive definiteness of the matrix
\[
A = \left(
\begin{array}{@{}ccccc@{}}
1 & -\e_{12} & -\e_{13} & \dots & -\e_{1n} \\
-\e_{21} & 1 & -\e_{23} & \dots & -\e_{2n} \\
-\e_{31} & -\e_{32} & 1 & \dots & -\e_{3n} \\
\vdots & \vdots & \vdots & \ddots & \vdots \\
-\e_{n1} & -\e_{n2} & -\e_{n3} & \dots & 1
\end{array}
\right),
\]
where $\e_{ij} = \cos \pi/m_{ij}$. This allows us to apply Theorem~\ref{th:main-distance-estimate}.

Let $v$ be $\e$-almost invariant unit vector in $\HH$ then by Theorem~\ref{th:main-distance-estimate}
the distance between $v$ and the space of $G$ invariant vectors in bounded by
\[
d_{\HH^G}(v)^2 \leq \left(\frac{\e}{2}\right)^2 \mathbf{1}^t A^{-1} \mathbf{1},
\]
where $\mathbf{1}$ is a column vector consisting of $n$ ones, because each generator $s_i$ moves
$v$ by $2 d_{\HH^{G_i}}(v) \leq \e$. Thus, if $\e < 2 (\mathbf{1}^t A^{-1} \mathbf{1})^{-1/2}$
there is a nontrivial invariant vector in $\HH$. This shows that
$\kc(G,S) \geq \e_0=2 (\mathbf{1}^t A^{-1} \mathbf{1})^{-1/2}$.
However, it is easy to see that there is an equality,
because the defining representation of $G$ contains an unit vector
which is
$\e_0$-almost invariant.

A similar argument can be used to obtain bounds for the spectral gap of the Laplacian:
The operator $\Id - \pi(s_i)$ is equal to two times the projection onto $\left(\HH^{G_i}\right)^\perp$.
Thus for a vector $v$ we have $(\Id - \pi(s_i)(v) = 2P_{\left(\HH^{G_i}\right)^\perp}(v)$, i.e.,
\[
\la \Delta_S(v),v \ra = \frac{1}{|S|} \sum \la 2 P_{\left(\HH^{G_i}\right)^\perp}(v) , v \ra
=
\frac{2}{|S|} \sum d_{\HH^{G_i}}(v)^2.
\]
If the vector $v$ has a trivial projection on the space of $G$ invariant vector we have
\[
\|v \|^2 = d_{\HH^G}(v)^2 \leq \mathbf{d}_v^t A^{-1}\mathbf{d}_v \leq
\lambda^{-1} \| \mathbf{d}_v \|^2 =
\lambda^{-1} \sum d_{\HH^{G_i}}(v)^2,
\]
where $\lambda$ is the smallest eigenvalue of the matrix $A$. Combining the above inequalities yields
\[
\la \Delta_S(v),v \ra \geq \frac{2}{n} \,\lambda \| v\|^2,
\]
i.e., the spectral gap of $\Delta_S$ is at least $\frac{2\lambda}{n}$. Again it is easy to see that there is
in equality since the smallest eigenvalue of $\Delta_S$ in the defining representation is
equal to $\frac{2\lambda}{n}$.
\end{proof}
\begin{example}
\label{ex:Coxeter_An}
If the Coxeter group $G$ is of type $A_n$, i.e., if $G\simeq \Sym(n+1)$
and $S$ consists of transpositions $(1,2), (2,3),\dots, (n,n+1)$.
In this case the matrix $A$ has the form
\[
\left(
\begin{array}{@{}rrrrr@{}}
1 & -\frac{1}{2} & 0 & \dots & 0 \\
-\frac{1}{2} & 1 & -\frac{1}{2} & \dots & 0 \\
0 & -\frac{1}{2} & 1 & \dots & 0 \\
\vdots\, & \vdots\, & \vdots\, & \ddots & \vdots\, \\
0 & 0 & 0 & \dots & 1
\end{array}
\right).
\]
A standard computation shows that this matrix has eigenvalues $\lambda_k=2\sin^2(\frac{k\pi}{2n+2})$
with eigenvectors
$v_k = (\dots \sin\frac{k \pi i}{n+1} \dots)^t$. Thus, the spectral gap %
of the Laplacian is
\[
\frac{2\lambda_1}{n} = \frac{4}{n} \sin^2 \frac{\pi}{2n+2} \sim \frac{\pi^2}{n^3},
\]
which implies that the relaxation time of the random walk on symmetric group is of order $n^3$.

The eigenvalues and the eigenvectors of $A$ can be used to compute the value %
of
$\mathbf{1}^t A^{-1} \mathbf{1}$, however it is easier%
\footnote{It is possible to bypass this computation by
constructing a vector in the defining representation, which at the same distance form each of the
fixed subspaces and use it to evaluate $\mathbf{1}^t A^{-1} \mathbf{1}$.}
to write down
explicitly the  matrix $A^{-1}$, and calculate that
$\mathbf{1}^t A^{-1} \mathbf{1}= (n^3+ 3n^2 + 2n)/6$, which implies that the Kazhdan constant of
symmetric group is
\[
\kc\left(\Sym(n+1),\{(1,2),(2,3),\dots,(n,n+1)\right) = \sqrt{\frac{24}{n^3+ 3n^2 + 2n}}.
\]
\end{example}
\begin{remark}
For any Coxeter group the smallest eigenvalue of the matrix $A$ is equal to
$2\sin^2 \frac{\pi}{2h}$, where $h$ is the Coxeter number of $G$.
This implies that the spectral gap of $\Delta_S$ is equal to
$
\displaystyle
\frac{2}{n}\left( 1 - \cos \frac{\pi}{h}\right)=\frac{4}{n} \sin^2\frac{\pi}{2h}
$.

The Kazhdan constant $\kc\left(G,S\right)$ is equal to
$2 \left(\mathbf{1}^t A^{-1} \mathbf{1}\right)^{-1/2}$.
In the simply laced case the number $M= \mathbf{1}^t A^{-1} \mathbf{1}$
is equal to the Dynkin index~\cite{Dyn}
of the canonical embedding of $\mathfrak{sl}_2$ in the simple Lie algebra
corresponding to the Coxeter group $G$ and by~\cite{Pan} is equal to $n h(h+1)/6$.
We do not know of any similar formula in the non-simply laced case.
However, it is not difficult to
compute the Kazhdan constant in each case
$\kc\left(G,S\right) =  2 M^{-1/2}$, where $M =\mathbf{1}^t A^{-1} \mathbf{1}$
is given in Table~\ref{table:Kazdan}.
\begin{table}
$$
\begin{array}{||c|c|c|c|c||}
\hline
\mbox{type} & \mbox{rank} & \mbox{Coxeter num.} & M  & \mbox{order of $M$}\vsp  \\
\hline
A_n & n & n+1 &
  n(n+1)(n+2)/6 & n^3/6 \\
B_n & n & 2n  &
  n\left(2n^2 + 3(\sqrt{2}-1)n + 4 -3 \sqrt{2}\right)/3 & 2n^3/3\\
D_n & n & 2(n-1) &
  n(n - 1)(2n - 1)/3 & 2n^3/3\\
E_6 & 6 & 12 &
 156 & \\
E_7 & 7 & 18 &
399 & \\
E_8 & 8 & 30 &
 1240 & \\
F_4 & 4 & 12 &
56 + 36 \sqrt{2} & \\
H_3 & 3 & 10 &
31 + 12 \sqrt{5} & \\
H_4 & 4 & 30 &
332 + 144 \sqrt{5} & \\
I_2(m) & 2 & m  &
2 (1-\cos \frac{\pi}{m})^{-1} & 4m^2/\pi^2\\
\hline
\end{array}
$$
\caption{Kazhdan constants for Coxeter groups.}
\label{table:Kazdan}
\end{table}
\end{remark}
\begin{remark}
Bacher and de la~Hapre~\cite{BdH} computed the Kazhdan constant
$\kc(G,S)$ when the Coxeter group $G$ is of type%
\footnote{They also considered the type $I_2(m)$,
 when $G$ is the dihedral group $D_m$.}
$A_n$, i.e., $G = \Sym(n+1)$.
Their proof uses
representation theory of the symmetric group and character estimates.
Bounds for the spectral gap of the Laplacian on the symmetric group ware obtained by
Diaconis and Shahshahani~\cite{DSh} and Diaconis and Saloff-Coste~\cite{DS}.

Bagno~\cite{Bag} extended the methods from~\cite{BdH}
to the case of Coxeter groups of types $B_n$ and $D_n$, but he was
not able to compute the exact value of the Kazhdan constants.
For the exceptional Coxeter groups
results of this type can be verified by long computation.
\end{remark}
\begin{remark}
It seems that Theorem~\ref{th:main-app-Coxeter}
can be generalized to finite complex reflection groups.
In order to do that one first needs prove that for any
unitary representation of any any rank 2 complex reflection
group the angle between the fixed subspaces of the two
generating psudo-reflections group is the same as the angle
in the defining representation. Since this is clearly the case
for the classical rank 2 complex reflection groups, one can
easily extends Theorem~\ref{th:main-app-Coxeter} the groups of type
$G_{m,p,n}$.
\end{remark}

\begin{remark}
One of the reasons why Theorem~\ref{th:main-app-Coxeter} gives exact bounds for the
Kazhdan constants and spectral gap is the existence of the defining representation, where
each ``generating'' subgroup fixes a hyperplane and the angles between these hyper planes
are lower bounds for the angles between the fixed subspaces of these subgroups in any
unitary representation of the group $G$. Theorem~\ref{th:main-app-RW-SO_n} is an other example,
where a similar configuration of subspaces allow us to compute the exact value of the spectral gap.
\end{remark}

\subsection{Property T for  Steinberg groups}
Theorem~\ref{th:main-distance-estimate} can be used to prove property T for groups $G$ which has
generating set $S$ consisting of several (finite) subgroups $S = \bigcup G_i$.
Before applying it, one only needs
to understand the representation theory of the subgroups of $G$ generated by any pair of $G_i$ and $G_j$.
One situation, where this method works nicely is the following:
\begin{theorem}
\label{th:main-T_for_St_groups}
The  Steinberg group $\St_n(\F_p \la t_1,\dots,t_k\ra)$ has property T,
provided that $n\geq 3$ and $p\geq 5$.
\end{theorem}
\begin{proof}
For an associative ring $R$ the  Steinberg group $\St_n(R)$ is generated by the elements $x_{ij}(r)$ where $1\leq i \not=j \leq n$
and $r\in R$ subject to the defining relations
\[
x_{ij}(r_1)x_{ij}(r_2) = x_{ij}(r_1+r_2)
\quad
[x_{ij}(r),x_{jk}(s)]=x_{ik}(rs)
\quad
[x_{ij}(r),x_{kl}(s)] =1.
\]

The group $G=\St_n(\F_p \la t_1,\dots,t_k\ra)$ contains $n$ finite subgroups $G_1,\dots, G_n$:
for $i=1,\dots,n-1$ the subgroup $G_i$ consists of the elements $x_{i,i+1}(a \cdot 1)$ for $a\in \F_p$
and the group $G_n$ consists of $x_{n1}(a_0 1 + a_1 t_1 + \dots +a_k t_k)$ for $a_i\in\F_p$.
An easy computation by induction shows that the subgroups $G_i$ generate the group $G$.

If $i,j < n$ and $|i - j | > 1$ then the subgroup groups $G_i$ and $G_j$ commute, therefore for
any unitary representation of $G$ their fixed subspaces are perpendicular, i.e.,
$\an{\HH^{G_i}}{\HH^{G_j}} = \pi /2$.
If $j=i+1 < n$ the $G_i$ and $G_j$ generate a Heisenberg group $H_p$ of order $p^3$.
It can be shown~\cite[Section 4]{EJ}
that in any representation of $H_p$ the angle between the fixed subspaces of any two non-central
cyclic subgroups of order $p$ is more $\cos^{-1}(p^{-1/2})$, thus
$\cos \an{\HH^{G_i}}{\HH^{G_{i+1}}} \leq p^{-1/2}$.

If one of the subgroups is $G_n$ the argument is almost the same---$G_n$
commutes with $G_2,\dots,G_{n-2}$ and the groups $\la G_1,G_n \ra$ and
$\la G_{n-1},G_n\ra$ are generalized Hisenberg groups. Thus, we have obtained bounds for all angles
$\an{\HH^{G_i}}{\HH^{G_j}}$ and the matrix $A$ has the form
\begin{equation}
\label{eq:affine_matrix}
\left(
\begin{array}{@{}cccccc@{}}
1 & -\e & 0 & \dots & 0 & -\e \\
-\e & 1 & -\e & \dots & 0 & 0 \\
0 & -\e & 1 & \dots & 0 & 0 \\
\vdots & \vdots & \vdots & \ddots & \vdots & \vdots \\
0 & 0 & 0 & \dots & 1 & -\e \\
-\e & 0 & 0 & \dots & -\e &1
\end{array}
\right),
\end{equation}
where $\e = p^{-1/2}$. If $p \geq 5$ this matrix is positive definite and its smallest eigenvalue is
equal to $1-2 \e$ with the corresponding eigenvector is $(1,1,\dots, 1)^t$.

By Theorem~\ref{th:main-distance-estimate} we have that for any $v\in \HH$ we have
\[
d_{\HH^G}(v)^2 \leq \frac{1}{1 - 2p^{-1/2}}\sum d_{\HH^{G_i}}(v)^2
\]
If the unit vector $v$ is $\epsilon$-almost invariant under the generating set $\bigcup G_i$ then
$d_{\HH^{G_i}}(v) \leq \epsilon/\sqrt{2}$ for each $i$,
because any unit vector in an unitary representation
of a group $H$ is moved by more than $\sqrt{2}$ by some element in $H$. %
Thus
\[
d_{\HH^G}(v)^2 \leq \frac{1}{1 - 2p^{-1/2}} \times \frac{n \epsilon^2}{2}.
\]
If $\epsilon < \sqrt{\frac{2(1-2p^{-1/2})}{n}}$ we have that $d_{\HH^G}(v) <1$, i.e., $\HH$ has
invariant vectors. This shows that
\[
\kc(G, \cup G_i) \geq \sqrt{\frac{2(1-2p^{-1/2})}{n}} \geq \sqrt{\frac{1}{5n}} > 0,
\]
in particular $G$ has Kazhdan property T, since $\cup G_i$ is finite.
\end{proof}

\begin{remark}
Theorem~\ref{th:main-T_for_St_groups} can be generalized to (hight rank)  Steinberg groups of other types
(over commutative rings) and the proof is essential the same---the groups $G_i$ for $i=1,\dots, n-1$
are part of the roots subgroups corresponding to the simple roots and $G_{n}$ is a subgroup of the root
subgroup corresponding to the largest negative root. In the simply laced case the
matrix $A$ is related to the Cartan matrix corresponding to the
extended Dynkin diagram---it will be positive definite, if $p\geq 5$
and the smallest eigenvalue will be again equal to $1 - 2p^{-1/2}$.
\end{remark}

\begin{remark}
One can use the method in the proof of Theorem~\ref{th:main-T_for_St_groups} to show that
the simply laced Kac-Moody groups over finite fields corresponding to $k$-regular graphs
have property T if $p \geq k^2$---there groups are generated by subgroups $G_i$,
indexed by the vertices of the graph. Each group $G_i$ is isomorphic to
$\SL_2(\F_p)$ and the group generated by $G_i$ and $G_j$ is either isomorphic to
$\SL_2(\F_p) \times \SL_2(\F_p)$ or $\SL_3(\F_p)$
depending whether the vertices are connected or not.
These conditions lead to bounds for the angles between the fixed spaces for $G_i$.
By Theorem~\ref{th:main-angle-bound} and Observation~\ref{ob:main_observation} the
positive definiteness of the matrix $A_\e$ implies that the resulting Kac-Moody group has property T.
However, the matrix $A_\e$ is positive define because its
smallest eigen value is equal to $1-kp^{-1/2} > 0$.
\end{remark}

\begin{proof}[Proof of Theorem~\ref{th:main-app-SL_n}]
The group $\SL_n(\F_p[t_1,\dots,t_k])$ is a quotient of the Steinberg group
$\St_n(\F_p\la t_1,\dots,t_k\ra)$ which has property T by Theorem~\ref{th:main-T_for_St_groups}.
However, property T is inherited by quotients, therefore $\SL_n(\F_p[t_1,\dots,t_k])$ also
has Kazhdan property T.
\end{proof}
\begin{remark}
This is not the first proof that the group $\St_n(\F_p\la t_1,\dots, t_n \ra)$ has property T.
The case $k=1$ is very old and goes back to Kazhdan~\cite{Kaz}---in this case $G$ has T because it is a
lattice in a high rank lie group over $\F_p((t^{-1}))$. In the case of commuting variables
the author and N. Nikolov~\cite{KN}
have shown that the group has property $\tau$, which is a weak form of property T. Also, in
the commutative
case Y. Shalom~\cite{Sh} proved that $G$ has T if $k \leq n-2$. This condition was replaced
by Vaserstein~\cite{Va} with $n\geq 3$. Recently Ershov and Jaikin~\cite{EJ}
extended these results by showing that
$\St_n(\Z\la t_1,\dots,t_k \ra)$ has property T for $n\geq 3$.

Essentially, the same proof as above valid in the case
$p > n^2 $ 
can also be found
in~\cite{EJ}. 
As mentioned in the introduction the aim of this paper is not
to prove new results but to explain the author's interpretation of the ideas
in~\cite{DJ} and~\cite{EJ}.
\end{remark}

\begin{remark}
Using results about relative property T one can also works with group $G_i$ which are not finite.
This approach yields that the  Steinberg groups $\St_n(R)$ has property T if $n\geq 3$, the ring $R$
is finitely generated and $R$ does not have $\F_q$ as a quotient, for $q \leq 4$.
It is possible to remove the condition that $R$ does not have small quotients,
however this requires significantly more complicated arguments instead of
Theorem~\ref{th:main-angle-bound}, see~\cite{EJ, EJK}.
\end{remark}



\subsection{Spectral gaps and mixing times of some random walks} 
The proofs of the following two theorems are very similar to the proofs of Theorems~\ref{th:main-app-Coxeter}
and Theorem~\ref{th:main-T_for_St_groups} and we will only sketch the main steps of the proofs.

\begin{theorem}
\label{th:main-app-RW-SL_n}
Let $G$ denote the group $\SL_n(\F_p)$ for $n\geq 3$ and $p\geq 5$
and let $G_{i,j}$ be the root subgroup $\Id + \F_p e_{ij}$ inside $G$.
The spectral gap of $\Delta_S$ is bounded by:
\begin{itemize}
\item[a)] $1/10n$ if $S = G_{1,2}\cup G_{2,3} \cup \dots \cup G_{n-1,n} \cup G_{n,1}$;
\item[b)] $1/10n$ if $S = \cup G_{i,j}$; 
\item[c)] $1/200n$ if $S = \{ \Id \pm e_{ij} \mid |i-j| \leq 1 (mod n) \}$. 
\end{itemize}
These bounds imply that the random walks of $\SL_n(\F_p)$ with
respect to the generating set described above have mixing time bounded by $C n^3 \log p$.
\end{theorem}
\begin{proof}
a) The proof is essentially the same as the one of Theorem~\ref{th:main-T_for_St_groups}.
Every pair of groups $G_i$ and $G_j$ either commute or generate a Heisenberg group.
The resulting matrix is the same as the one in (\ref{eq:affine_matrix}) and its smallest eigenvalue
is equal to $1-2p^{-1/2}$. Now using that
\[
\left\langle \frac{1}{|G_i|} \sum_{g\in G_i} \left( \Id -\pi(g)\right)v ,v \right\rangle = d(v,\HH^{G_i})^2
\]
we obtain that the spectral gap of $\Delta$ is bounded by
\[
\Delta_S > \frac{1}{n} \left(1-2p^{-1/2} \right) > \frac{1}{10n}.
\]

b) Follows from part a) by writing the decomposing the complete graph as union of $n$-cycles and
observing that the Laplacian is the average of the Laplacians corresponding the union of
the root subgroups in each cycle.

c) First, one use the Sleberg's theorem to bound  the spectral gap of the Lapalcian on
$\SL_2(\F_p)$ with respect to the generating set consisting of $\Id \pm e_{12}$ and $\Id \pm e_{12}$.
This bound implies that 
if a vector $v$ in any unitary representation of $\SL_n(\F_p)$
is $\e$-almost invariant with respect to $S$ then it is $20\e$-almost invariant with respect to
$\bigcup_{|i-j| =1} G_{ij}$, which combined with part a) completes the proof.
\end{proof}
\begin{remark}
Theorem~\ref{th:main-app-RW-SL_n} implies that the Kazhdan constant of $\SL_n(\F_p)$ with respect
to the generating sets $S = \bigcup G_{i,j}$ is of the bounded below by a function of order $n^{-1/2}$. A bound of this type was found in~\cite{KSLn}, however this argument gives slightly better constant. It is easy to construct
representations of $\SL_n(\F_p)$ with $n^{-1/2}$-almost invariant vectors which shows that there is in upper
bound for the Kazhdan constant of the same order.
\end{remark}

\medskip

Our final example is slightly different, because in involves a compact Lie group.
Let $G$ denote the group $\SO(n)$ (with its standard action on $\R^n$).
This group contains the subgroups
$G_{ij}\simeq \SO(2)$ consisting of rotatoins in the coordionate plane $\mathrm{span}\{e_i,e_j\}$.
Kac~\cite{Kac} studied
the random walk on $\SO(n)$ with respect to $\bigcup G_{ij}$.
Maslen~\cite{Mas} computed the spectral gap in the case
$
S=\left\{ G_{ij} \mid \vsp {1\leq i,j \leq n} \right\}
$
and
Diaconis and Saloff-Coste~\cite{DS2} obtained bounds in the case
$
S=\left\{ G_{i,i+1} \mid \vsp {1\leq i <n} \right\}
$.
Since the group $G$ is not finite one need to slightly modify the definition of the Laplacian
$\Delta_S$---instead of averaging over the generating set one uses an integral with respect to some
measure $\mu$, in this case $S$ is a union of circles and $\mu$ is just the average of the
uniform measures on each circle.

\begin{theorem}
\label{th:main-app-RW-SO_n}
The spectral gap of $\Delta_S$ 
is
\begin{itemize}
  \item[a)]
     is equal to
     $\delta$
     if $S = \bigcup_{1\leq i \leq n} G_{i,i+1}$;
  \item[b)]
     bounded below by
     $\delta$
     if $S = \bigcup_{1\leq i,j \leq n+1} G_{ij}$,
\end{itemize}
where $\displaystyle \delta = \frac{2}{n} \sin^2\left( \frac{\pi}{2n+2}\right) \sim \frac{\pi^2}{2n^3}$.
\end{theorem}
\begin{remark}
Part a) improves the bound found in~\cite{DS2} by a constant factor.
On the other side the bound in part b) is significantly weaker then
exact value of the gap
$\displaystyle
\frac{n+3}{2n(n+1)}\sim\frac{1}{2n}$ see~\cite[Teorem 2.1]{Mas}
\end{remark}
\begin{proof}
The proof of part a) is similar to Theorem~\ref{th:main-app-Coxeter}
and Example~\ref{ex:Coxeter_An}.  The measure
$\mu_{G_{ij}}$ on $G_{ij}$ is uniform, therefore 
\[
\int_{\mu_{G_{ij}}} \left(\vsp \Id - \pi(g)\right)(v)\, d\mu =
P_{\left(\HH^{G_{ij}}\right)^\perp}(v).
\]
Therefore, if we denote $G_i= G_{i,i+1}$, we have
\[
\la \Delta_S(v),v \ra = \frac{1}{n} \sum \la P_{\left(\HH^{G_i}\right)^\perp}(v) , v \ra
=
\frac{1}{n} \sum d_{\HH^{G_{i}}}(v)^2.
\]
The group $G_i$ and $G_j$ commute if $|i-j|>1$ therefore
$\cos\an{\HH^{G_i}}{\HH^{G_j}} \leq 0$ (the representation $\HH$ might not have
any invariant vectors under $G_i$, in which case the angle will be equal to $\pi$).
If $j=i+1$ the groups $G_i$ and $G_j$ generate a group isomorphic to $\SO(3)$.
Using the representation theory of $\SO(3)$ one can show~\cite[Lemma 3.2]{Mas} that
$\cos\an{\HH^{G_i}}{\HH^{G_j}} \leq 1/2$.

Thus, the matrix $A$ is the same one in Example~\ref{ex:Coxeter_An},
and its smallest eigenvalue is equal to $\lambda=2\sin^2\left(\frac{\pi}{2n+2}\right)$,
which implies that the spectral gap of $\Delta_S$ is bounded bellow by
$$
\frac{2}{n} \sin^2\left( \frac{\pi}{2n+2}\right) \sim \frac{\pi^2}{2n^3}.
$$

Actually there is an equality, because spectral gap of the Laplacian
on the representation of $\SO(n+1)$ on the space of harmonic
homogeneous polynomials of degree $2$, is exactly equal to
$\delta$---this representation contains a subspace $V$ of dimension $n$ which contains
$n$ hyper planes $\HH^{G_i}$ and the angles between $\HH^{G_i}$ and $\HH^{G_j}$
is either $\pi/2$ or $\pi/3$.

Part b) follows immediately from part a) by writing the generating set as union of several
generating sets for part a).
\end{proof}


\section*{Acknowledgments}
The author thanks Mikhail Ershov and Andrei Jaikin-Zapirain for the useful discussions
about the paper, especially for explaining the ideas behind~\cite{DJ} and~\cite{EJ}.
The main result started as a conjecture by the author and M. Ershov formulated during Ershov's
visit to Cornell University in December 2008.
The paper was written during the author's visit to the Institute for Advanced Study in the Fall 2009.


\bigskip
\noindent
\textsc{%
Department of Mathematics,
Cornell University, \\
Ithaca, NY 14853, USA
\\
and
\\
School of Mathematics,
Institute of Advanced Study,\\
Princeton, NJ 08540, USA
}

\noindent
email: {\tt kassabov@math.cornell.edu}

\end{document}